\documentclass[12pt]{article}
\usepackage{amssymb, amsmath, latexsym, a4}
\usepackage[latin1]{inputenc}
\usepackage[all]{xy}
\usepackage[usenames]{color}
\begin{document}

\newcommand{\rdg}{\hfill $\Box $}

\newtheorem{De}{Definition}[section]
\newtheorem{Th}[De]{Theorem}
\newtheorem{Pro}[De]{Proposition}
\newtheorem{Co}[De]{Corollary}
\newtheorem{Rem}[De]{Remark}
\newtheorem{Ex}[De]{Example}
\newtheorem{Exo}[De]{Exercises}
\newtheorem{Le}[De]{Lemma}
\newcommand{\tp}{\otimes}
\newcommand{\N}{\mathbb{N}}
\newcommand{\Z}{\mathbb{Z}}
\newcommand{\K}{\mathbb{K}}

\newcommand{\cok}{{\sf Coker}}
\newcommand{\Hom}{{\sf Hom}}
\newcommand{\im}{{\sf Im}}
\newcommand{\ext}{{\sf Ext}}
\newcommand{\id}{{\sf id}}
\newcommand{\Ker}{{\sf Ker}}

\newcommand{\ele}{\cal L} \newcommand{\as}{\cal A} \newcommand{\ka}{\cal K}\newcommand{\eme}{\cal M} \newcommand{\pe}{\cal P}
\newcommand{\he}{\cal H} \newcommand{\ene}{\cal N} \newcommand{\qu}{\cal Q} \newcommand{\be}{\cal B}

\newcommand{\pn}{\par \noindent}
\newcommand{\pbn}{\par \bigskip \noindent}
\bigskip\bigskip

\centerline{\large \bf Universal $\alpha$-central extensions of Hom-Leibniz $n$-algebras}
\bigskip

\centerline{\bf J. M. Casas$^{(1)}$ and N. Pacheco Rego$^{(2)}$}

\bigskip \bigskip
\centerline{$^{(1)}$Departamento de  Matemática Aplicada I, Universidade de Vigo,}
\centerline{E. E. Forestal, 36005 Pontevedra, Spain}
\centerline{e-mail address: \tt jmcasas@uvigo.es}
\bigskip

\centerline{$^{(2)}$IPCA, Dpto. de Ciências, Campus do IPCA,
 Lugar do Aldão}
\centerline{4750-810 Vila Frescainha, S. Martinho, Barcelos,
 Portugal}
\centerline{e-mail address: \tt natarego@gmail.com}

\bigskip \bigskip \bigskip

\par
{\bf Abstract}
 We construct homology with trivial coefficients of Hom-Leibniz $n$-algebras.  We introduce and characterize universal ($\alpha$)-central extensions of Hom-Leibniz $n$-algebras. In particular, we show their interplay with the zeroth and first homology with trivial coefficients.
  When $n=2$  we recover the corresponding results on universal central extensions of Hom-Leibniz algebras. The notion of non-abelian tensor product of Hom-Leibniz $n$-algebras is introduced and we establish its relationship with  the universal central extensions. We develop a generalization of the concept and properties of unicentral Leibniz algebras to the setting of Hom-Leibniz $n$-algebras.

\bigskip \bigskip


\bigskip \bigskip
{\it A. M. S. Subject Class. (2000):} 17A30,  17B55, 18G60

\section{Introduction}
Algebras endowed with an $n$-ary operation play important roles, among others,  in Lie and Jordan theories, geometry, analysis,  physics and biology.
For instance, this kind of structures were considered to analyze  DNA recombination \cite{Sv}. Leibniz $n$-algebras and its corresponding skew-symmetric version, named as Lie $n$-algebras or Filippov algebras, arose in the setting of  Nambu mechanics \cite{Na}, a generalization of the Hamiltonian mechanics.
The particular case $n=3$ has found applications in  string theory and M-branes \cite{BL, PS} and in the M-theory generalization of the Nahm's equation proposed by Basu and Harvey \cite{BH}. It can also be used to construct solutions of the Yang-Baxter equation \cite{OS}, which first appeared in statistical mechanics \cite{BRJ}.

Deformations of algebras structures by means of endomorphisms give rise to Hom-algebra structures. They are motivated by discrete and deformed vector fields and differential calculus. Part of the reason to study Hom-algebras is its relation with the $q$-deformations of the Witt and the Virasoro algebras (see \cite{HLS}).

In this way, deformations of algebras of Lie type were considered, among others,  in \cite{HLS, MS, MS2, Yau, Yau1}. Deformations of algebras of Leibniz type were considered, among others, in \cite{ChS, GI, CIP1,  Is2, MS}. The generalizations of $n$-ary algebra structures, such as Hom-Leibniz $n$-algebras (or n-ary Hom-Nambu) and Hom-Lie  $n$-algebras (or n-ary Hom-Nambu-Lie), have been introduced in \cite{AMS} by Ataguema, Makhlouf, and Silvestrov. In these Hom-type algebras, the n-ary Nambu identity is deformed using $n-1$ linear maps, called the twisting maps, given rise to the fundamental identity (n-ary Hom-Nambu identity) (see Definition \ref{HomLeibn}). When these twisting maps are all equal to the identity map, one recovers Leibniz $n$-algebras (n-ary Nambu) and Lie $n$-algebras (Nambu-Lie algebras).

The topic of central extensions of algebraic structures is also present in many applications to Physics. For instance, the Witt algebra and its  one-dimensional universal central extension, the Virasoro algebra, often appear in problems with conformal symmetry in the setting of string theory \cite{Gan}.

Recently in \cite{CVdL} was noticed an important fact concerning universal central extensions in the setting of semi-abelian categories, the so called {\bf UCE} condition, namely the composition of two central extensions is central. We show in this paper that the category of Hom-Leibniz $n$-algebras doesn't satisfy  {\bf UCE}  condition (see Example \ref{uce condition}). From this fact, our aim in this article is to introduce and characterize universal $\alpha$-central extensions of Hom-Leibniz $n$-algebras. In case $n=2$ we recover the corresponding results on universal $\alpha$-central extensions of Hom-Leibniz algebras in \cite{CIP1, CKP1}. Moreover, in case $\alpha = \id$ we recover results on universal central extensions of Leibniz $n$-algebras in \cite{Ca2}. In case $n=2$ and  $\alpha = \id$ we recover results from \cite{CC}.

The article is organized as follows: in  section \ref{preliminaries}  we introduce the necessary basic concepts on  Hom-Leibniz $n$-algebras and construct the homology with trivial coefficients of Hom-Leibniz $n$-algebras. Bearing in mind \cite{Ca1}, we endow the underlying vector space to a Hom-Leibniz $n$-algebra ${\ele}$ with a structure of $({\cal D}_{n-1}\left( {\cal L} \right) = {\ele}^{\otimes n-1}, \alpha')$-symmetric Hom-co-representation as Hom-Leibniz algebras and define the homology  with trivial coefficients of ${\ele}$ as the Hom-Leibniz homology $HL_{\ast}^{\alpha}({\cal D}_{n-1}\left( {\ele} \right), {\ele})$.

In section \ref{uce} we present our main results on universal central extensions. Based on the investigation initiated in \cite{CIP1}, we generalize the concepts of ($\alpha$)-central extension, universal ($\alpha$)-central extension and perfection to the framework of Hom-Leibniz $n$-algebras. We also extend the corresponding characterizations of universal ($\alpha$)-central extensions.
Since Hom-Leibniz $n$-algebras category doesn't satisfy  {\bf UCE}  condition,  characterizations are divided between universal central and universal $\alpha$-central (see Theorem \ref{teorema}).

In section \ref{non-abelian tensor} we introduce the concept of non-abelian tensor product of Hom-Leibniz $n$-algebras that generalizes the non-abelian tensor product of Leibniz ($n$)-algebras in \cite{Ca2, CKP1}, and we establish its relationship with  the universal central extension.

The final section is devoted to develop a generalization of the concept and properties of unicentral Leibniz algebras in \cite{CC} to the setting of Hom-Leibniz $n$-algebras. As a first step we show that the classical result: perfect Leibniz algebras are unicentral, doesn't hold in the framework of Hom-Leibniz $n$-algebras (see Example \ref{unicentral}) and requires an additional condition (see Proposition \ref{prop 1}). The main result in this section establishes that for two perfect Hom-Leibniz $n$-algebras, $({\ele}, \widetilde{\alpha}_{\ele})$ and $ ({\ele}', \widetilde{\alpha}_{\ele'})$  with both $\alpha_{\ele}, \alpha_{\ele'}$ injective and such that  $\left( \frak{uce}({\ele}), \widetilde{\alpha}_{\frak{uce}({\ele})} \right)$, and $\left( \frak{uce}({\ele}'), \widetilde{\alpha}_{\frak{uce}({\ele}')} \right)$  satisfy condition (\ref{condition}) (see below), then the following statements hold:
 \begin{enumerate}
 \item[a)] If  $\left( \frak{uce}({\ele}), \widetilde{\alpha}_{\frak{uce}({\ele})} \right) \cong \left( \frak{uce}({\ele}'), \widetilde{\alpha}_{\frak{uce}({\ele}')} \right)$, then $\frac{\alpha_{\ele}({\ele})}{Z( \alpha_{\ele}({\ele}))} \cong \frac{\alpha_{\ele'}({\ele'})}{Z( \alpha_{\ele'}({\ele'}))}$.
 \item[b)] If $\frac{\alpha_{\ele}({\ele})}{Z( \alpha_{\ele}({\ele}))} \cong \frac{\alpha_{\ele'}({\ele'})}{Z( \alpha_{\ele'}({\ele'}))}$, then $\left( \frak{uce}(\alpha_{\ele}({\ele})), \widetilde{\alpha}_{\frak{uce}({\ele})\mid} \right) \cong \left( \frak{uce}(\alpha_{\ele'}({\ele'})), \widetilde{\alpha}_{\frak{uce}({\ele}')\mid} \right)$.
 \end{enumerate}


\section{Preliminaries on Hom-Leibniz $n$-algebras} \label{preliminaries}

In this section we introduce necessary material on  Hom-Leibniz $n$-algebras, also called $n$-ary Hom-Nambu algebras in \cite{AMM1, AMS, Yau2} or $n$-ary Hom-Nambu-Lie algebras in \cite{AAS}.

\subsection{Basic definitions}

\begin{De}\label{HomLeibn}
A Hom-Leibniz $n$-algebra is triple $\left(  {\cal L},\left[-,\dots,-\right], \widetilde{\alpha }\right)$ consisting of a $\mathbb{K}$-vector space $\cal L$ equipped with an $n$-linear map $\left[ -,\dots,-\right] : {\cal L}^{\times n}\longrightarrow {\cal L}$ and a family $\widetilde{\alpha }=\left( \alpha _{i}\right), {1\leq i\leq n-1}$ of linear maps $\alpha _{i}: {\cal L} \longrightarrow {\cal L}$, satisfying the following fundamental identity:
\begin{equation}\label{def}
\begin{array}{c}
[[ x_{1},x_{2},\dots,x_{n}],\alpha _{1}(y_{1}),\alpha _{2}(y_{2}),\dots,\alpha _{n-1}(y_{n-1})]
=\\
\underset{i=1}{\overset{n}{\sum }}[\alpha _{1}(x_{1}),\dots,\alpha _{i-1}( x_{i-1}),[x_{i},y_{1},y_{2,}, \dots,y_{n-1}] ,\alpha_{i}(x_{i+1}),\dots,\alpha _{n-1}(x_{n})]
\end{array}
\end{equation}
for all $\left( x_{1},\dots,x_{n}\right) \in {\cal L}^{\times n}, y=\left(y_{1},\dots,y_{n-1}\right) \in {\cal L}^{\times (n-1)}$.
\end{De}

The linear maps $\alpha_1, \dots, \alpha_{n - 1}$ are called the twisting maps of the Hom-Leibniz $n$-algebra.
When the $n$-ary bracket is skew-symmetric, i.e. $[x_{\sigma(1)}, \dots, x_{\sigma(n)}] =$ $(-1)^{\epsilon(\sigma)}[x_1, \dots, x_n], \sigma \in S_n$, then the structure is called Hom-Lie $n$-algebra (or $n$-ary Hom-Nambu algebras in \cite{AMM, AAS}, or $n$-Hom-Lie algebras \cite{KMS}).
\bigskip

Let $x=\left( x_{1},\dots,x_{n}\right) \in {\cal L}^{\times n}$, $y=\left( y_{1},\dots,y_{n-1}\right) \in {\cal L}^{\times (n-1)}$, $\widetilde{\alpha }\left( y\right) =\left( \alpha _{1}\left( y_{1}\right),\dots, \right.$ $\left. \alpha_{n-1}\left( y_{n-1}\right) \right) \in {\cal L}^{\times (n-1)}$ and define the adjoint representation as the linear map $ad_{y}: {\cal L} \longrightarrow {\cal L}$, such that $ad_{ y}\left( x\right) =\left[ x, y_{1},\dots,y_{n-1} \right]$, for all $y \in {\cal L}$. Then  identity (\ref{def}) may be written as follows:
$$ad_{\widetilde{\alpha }\left( y\right)}[x_1, \dots, x_n] =\underset{i=1}{\overset{n}{\sum }}\left[ \alpha _{1}\left( x_{1}\right),\dots,\alpha_{i-1}\left( x_{i-1}\right) ,ad_{y}\left( x_{i} \right) ,\alpha _{i}\left( x_{i+1}\right),\dots,\alpha _{n-1}\left( x_{n}\right) \right]$$

\begin{De} \cite{AMM}
A Hom-Leibniz $n$-algebra $({\cal L},[-,\dots,-],\widetilde{\alpha })$  is said to be multiplicative if the linear maps in the family  $\widetilde{\alpha }=\left( \alpha _{i}\right) _{1\leq i\leq n-1}$ are of the form $\alpha_{1}= \dots = \alpha_{n-1}=\alpha$, and they preserve the bracket, that is,  $\alpha[x_1,\dots,x_n] = [\alpha(x_1),\dots,\alpha(x_n)]$, for all $(x_1,\dots,x_{n}) \in {\cal L}^{\times n}$.
\end{De}

\begin{De}\label{homo} \cite{AMM}
A homomorphism between two Hom-Leibniz $n$-algebras $\left(  {\cal L}, \right.$ $\left. \left[-,\dots,-\right], \widetilde{\alpha }\right)$ and $\left(  {\cal L}',\left[-,\dots ,-\right]', \widetilde{\alpha'}\right)$ where $\widetilde{\alpha }=\left( \alpha _{i}\right)$ and $\widetilde{\alpha}'=\left( \alpha_{i}'\right), 1\leq i\leq n-1$, is a linear map $f : {\cal L} \to {\cal L}'$ such that:
\begin{enumerate}
\item[(a)] $f([x_1, \dots , x_n]) = [f (x_1), \dots , f (x_n)]'$;
\item[(b)] $f \circ \alpha_i = \alpha_i'\circ f, i = 1, \dots , n - 1$
\end{enumerate}
for all $x_1,...,x_n \in {\ele}$.
\end{De}
\medskip

We denote by ${\sf _{n}HomLeib}$ the category of Hom-Leibniz $n$-algebras. In case  $n=2$, identity (\ref{def}) is the Hom-Leibniz identity (2.1) in \cite{CIP1}, so  Hom-Leibniz $2$-algebras are exactly  Hom-Leibniz algebras and we use the notation {\sf HomLeib} instead of ${\sf {_2}HomLeib}$.

\begin{Ex}\label{ejemplo 1} \
\begin{enumerate}
\item[(a)] When the maps $\left( \alpha _{i}\right) _{1\leq i\leq n-1}$ in Definition \ref{HomLeibn} are all of them the identity maps, then one recovers the definition of  Leibniz $n$-algebra \cite{CLP}. Hence Hom-Leibniz $n$-algebras include Leibniz $n$-algebras as a full subcategory, thereby motivating the name "Hom-Leibniz $n$-algebras" as a deformation of Leibniz $n$-algebras twisted by homomorphisms. Moreover it  is a multiplicative Hom-Leibniz $n$-algebra.

    \item[(b)] Hom-Lie $n$-algebras  are Hom-Leibniz $n$-algebras whose bracket satisfies the condition $[x_1,\dots,x_i,x_{i+1},\dots,x_n]=0$ as soon as $x_i=x_{i+1}$ for $1\leq i\leq n-1$. So the category ${\sf {_n}HomLie}$ of  Hom-Lie $n$-algebras can be considered as a full subcategory of ${\sf {_n}HomLeib}$. For any multiplicative Hom-Leibniz $n$-algebra $({\cal L},[-,\dots,-],\widetilde {\alpha})$ there is associated the Hom-Lie $n$-algebra $({\cal L}_{\rm Lie},[-,\dots,-],$ $\overline{\widetilde{\alpha}})$, where ${\cal L}_{\rm Lie} = {\cal L}/{\cal L}^{\rm ann}$, the bracket is the canonical bracket induced on the quotient and  $\overline{\widetilde{\alpha}}$ is the homomorphism naturally induced by $\widetilde{\alpha}$. Here ${\cal L}^{\rm ann} = \langle \{[x_1,\dots,x_i, x_{i+1}, \dots, x_n], \text{as soon as}\  x_i = x_{i+1}, 1 \leq i \leq n-1, x_j \in {\cal L}, j = 1, \dots, n \} \rangle$.

\item[(c)] Any Hom-vector space $V$ together with the trivial $n$-ary bracket $[-,-,\dots,-]$ (i.e. $[x_1, x_2, \dots, x_n]= 0$ for all $x_i \in V, 1 \leq i \leq n$) and any collection of linear maps $\widetilde{\alpha}_V = (\alpha_i:V \to V)_{1 \leq i \leq n-1}$, is a Hom-Leibniz $n$-algebra, called abelian Hom-Leibniz $n$-algebra.

\item[(d)] Hom-Lie triple systems \cite{ANI, Yau2} are  Hom-Leibniz $3$-algebras ${\cal L}$ satisfying the following properties:
\begin{enumerate}
\item[$\bullet$]  $[x,y,z]=-[y,x,z]$,
\item[$\bullet$]  $[x,y,z]+[y,z,x]+[z,x,y]=0$,
\end{enumerate}
for all $x, y, z \in {\cal L}$.

\item[(e)] 1-dimensional Hom-Leibniz $n$-algebras over a field $\mathbb{K}$, whose characteristic in not a factor of $n-1$, are abelian  Hom-Leibniz $n$-algebras or Hom-Leibniz $n$-algebras with any bracket and the collection $\widetilde{\alpha} = (\alpha_i)_{1 \leq i \leq n-1}$ contains at least one trivial map $\alpha_i, 1 \leq i \leq n-1$.

\end{enumerate}
\end{Ex}

In the sequel we refer to Hom-Leibniz $n$-algebras as multiplicative Hom-Leibniz $n$-algebras and we shall use the shortened notation
$({\cal L},\widetilde{\alpha }_{\cal L})$ when there is not confusion with the bracket operation.

\begin{De}
Let $({\cal L},[-,\dots,-],\widetilde{\alpha }_{\ele})$ be a Hom-Leibniz $n$-algebra. A  Hom-Leibniz $n$-subalgebra $({\cal H}, \widetilde{\alpha }_{\cal H})$ is a linear subspace ${\cal H}$ of ${\cal L}$, which is closed for the bracket and invariant by $\widetilde{\alpha }_{\cal L}$, that is,
\begin{enumerate}
\item [a)] $[x_1,\dots,x_n] \in {\cal H},$ for all $x_1,\dots,x_n \in {\cal H}$,
\item [b)] $\alpha_{\cal H}(x) \in {\cal H}$, for all $x  \in {\cal H}$ ($\alpha_{\cal H} = \alpha_{{\cal L} \mid}$).
\end{enumerate}

A  Hom-Leibniz $n$-subalgebra $({\cal H},\widetilde{\alpha}_{\cal H})$ of $({\cal L},\widetilde{\alpha}_{\cal L})$ is said to be an $n$-sided Hom-ideal if $[x_1,x_2,\dots,x_n] \in {\cal H}$ as soon as $x_i \in {\cal H}$ and $x_1,\dots,x_{i-1},x_{i+1},\dots,x_n \in {\cal L}$, for all $i=1,2,\dots,n$.

If $({\cal H},\widetilde{\alpha }_{\cal H})$ is an $n$-sided Hom-ideal of $({\cal L},\widetilde{\alpha }_{\cal L})$, then the quotient $\left( {\cal L}/{\cal H}, \widetilde{\alpha}_{{\cal L}/{\cal H}} \right)$ naturally inherits a structure of Hom-Leibniz $n$-algebra, which is said to be the quotient Hom-Leibniz $n$-algebra.
\end{De}

\begin{De}
Let $({\eme}_i, \alpha_{\ele \mid}), 1 \leq i \leq n,$ be subalgebras of a Hom-Leibniz $n$-algebra $({\cal L},\widetilde{\alpha}_{\cal L})$. We call commutator subspace corresponding to the subalgebras ${\eme}_i, 1 \leq i \leq n$, to the vector subspace of ${\ele}$  $$\left[{\eme}_1, \dots,{\eme}_n \right]= \langle \left \{ [x_1, \dots, x_n], x_i \in {\eme}_{\sigma(i)}, 1 \leq i \leq n, \sigma \in S_n \right \} \rangle$$
\end{De}

\begin{De}
Let $({\cal L},\widetilde{\alpha}_{\cal L})$ be a Hom-Leibniz $n$-algebra. The subspace
\[
\begin{array}{lcl} Z({\cal L}) &=& \{ x \in {\cal L} \mid [x_1,\dots,x_{i-1},x,x_{i+1},\dots,x_n] =0, \\
&&  \forall x_j \in {\cal L},j \in \{1,\dots, \widehat{i}, \dots,n\}, i \in \{1, \dots, n\} \}
 \end{array} \]
 is said to be the center of $({\cal L},\widetilde{\alpha}_{\cal L})$.

When the endomorphism $\alpha: {\cal L} \to {\cal L}$ is surjective, then  $Z({\cal L})$ is an $n$-sided Hom-ideal of ${\cal L}$.
\end{De}

\begin{Pro} \cite[Theorem 4.8 (2)]{ZCh}
Let $({\cal L},\widetilde{\alpha}_{\cal L})$  be a Hom-Leibniz $(n+1)$-algebra. Then $\left({\cal D}_{n}({\cal L})={\cal L}^{\otimes n}, [-,-], \alpha' \right)$ is a Hom-Leibniz algebra with respect to the bracket
$$\left[ a_{1}\otimes \dots \otimes a_{n},b_{1}\otimes \dots \otimes b_{n}\right] :=\underset{i=1}{\overset{n}{\sum }}\alpha\left( a_{1}\right) \otimes \dots \otimes \left[ a_{i},b_{1}, \dots,b_{n}\right] \otimes \dots \otimes \alpha\left( a_{n}\right)$$
and endomorphism $\alpha'= {\cal D}_{n}({\cal L}) \to {\cal D}_{n}({\cal L})$ given by
$$\alpha'(a_{1}\otimes ...\otimes a_{n})=\alpha(a_1)\otimes ...\otimes \alpha(a_n).$$
\end{Pro}


\subsection{Homology with trivial coefficients of Hom-Leibniz $n$-algebras}

Let $({\cal L}, \widetilde{\alpha}_{\cal L})$ be a Hom-Leibniz $n$-algebra, then ${\cal L}$ (as a $\mathbb{K}$-vector space) is endowed with a symmetric Hom-co-representation structure \cite[Definition 3.1]{CIP1} over $({\cal D}_{n-1}\left( {\cal L} \right) ={\cal L}^{\otimes \left( n-1\right)}, \alpha')$ as Hom-Leibniz
algebras with respect to the following actions
\[
\begin{array}{llcl}
& \left[ -,-\right] &:& {\cal L} \times {\cal D}_{n-1}\left( {\cal L} \right) \longrightarrow {\cal L}\\
& \left[ -,-\right] &:& {\cal D}_{n-1}\left( {\cal L} \right) \times {\cal L} \longrightarrow {\cal L}
\end{array}
\]
\noindent given by
\begin{equation} \label{action}
\begin{array}{llcr}
& \left[ l,l_{1}\otimes \dots \otimes l_{n-1}\right] &:=& \left[ l,l_{1}, \dots, l_{n-1}\right]\\
& \left[ l_{1}\otimes \dots \otimes l_{n-1},l\right] &:=&-\left[ l,l_{1}, \dots, l_{n-1}\right]
\end{array}
\end{equation}

\noindent and endomorphism $\alpha : {\cal L} \to {\cal L}$ such that $\widetilde{\alpha}_{\ele} = (\alpha_i), \alpha_i= \alpha, 1 \leq i \leq n-1.$

Now we construct a chain complex for Hom-Leibniz $n$-algebras in order to compute its homology with trivial coefficients. Firstly
we recall this complex for Hom-Leibniz homology \cite {CIP1}.
Let $(L,[-,-],\alpha_L)$ be a Hom-Leibniz algebra and $(M,\alpha_M)$ be a Hom-co-representation over $(L,[-,-],\alpha_L)$. The Hom-Leibniz complex $(CL_{\ast}^{\alpha}(L,M), d_{\ast})$ is given by setting  $CL_n^{\alpha}(L,M) :=M \otimes L^{\otimes n}, n \geq 0$, and by differentials the $\mathbb{K}$-linear maps $d_n: CL_n^{\alpha}(L,M) \to CL_{n-1}^{\alpha}(L,M)$ defined by
$$d_n(m \otimes x_1 \otimes \dots \otimes x_n)= [m,x_1] \otimes  \alpha_L(x_2) \otimes \dots \otimes \alpha_L(x_n) +$$
$$\sum_{i=2}^n (-1)^i [ x_i,m] \otimes \alpha_L(x_1) \otimes \dots \otimes  \widehat{\alpha_L(x_i)} \otimes \dots \otimes  \alpha_L(x_n) +$$
 $$\sum _{1 \leq i < j \leq n} (-1)^{j+1} \alpha_M(m) \otimes \alpha_L(x_1)\otimes  \dots \otimes  \alpha_L(x_{i-1}) \otimes [x_i,x_j] \otimes \dots \otimes \widehat{\alpha_L(x_j)} \otimes \dots \otimes \alpha_L(x_n).$$

 The homology of the chain complex $(CL_{\ast}^{\alpha}(L,M), d_{\ast})$ is called homology of the Hom-Leibniz algebra  $(L,[-,-],\alpha_L)$ with coefficients in the Hom-co-representation $(M,\alpha_M)$ \cite{CIP1} and is denoted by $HL_{\ast}^{\alpha}(L,M) :=H_{\ast}(CL_{\ast}^{\alpha}(L,M), d_{\ast})$.

 In order to construct the  chain complex $({_n}CL_{\ast}^{\alpha}(\cal L), \delta_{\star})$ which allows  the computation of the homology with trivial coefficients of a Hom-Leibniz $n$-algebra $({\cal L}, \widetilde{\alpha}_{\cal L})$, we only need to have in mind that (\ref{action}) endows  $({\cal L}, \alpha)$ with a Hom-co-representation structure over  $({\cal D}_{n-1}\left( {\cal L} \right), \alpha')$, so it makes sense the construction of its Hom-Leibniz complex, hence  we define
$${_n}CL_{\ast}^{\alpha}({\cal L}):=CL_{\ast}^{\alpha}({\cal D}_{n-1}\left( {\cal L} \right), {\cal L})$$
thus, by definition, the homology with trivial coefficients for the Hom-Leibniz $n$-algebra $({\cal L}, \widetilde{\alpha}_{\cal L})$ is
$${_n}HL_{\ast}^{\alpha}({\cal L}, \mathbb{K}):=HL_{\ast}^{\alpha}({\cal D}_{n-1}\left( {\cal L} \right), {\cal L})$$
and we will use the short notation  ${_n}HL_{\ast}^{\alpha}({\cal L})$ instead of ${_n}HL_{\ast}^{\alpha}({\cal L}, \mathbb{K})$.

 In particular, we have
$${_n}HL_{0}^{\alpha}({\cal L}) = HL_0^{\alpha}\left( {\cal D}_{n-1}\left( {\cal L} \right), {\cal L} \right) = \cok(d_1 : {\cal L}^{\otimes n} \to {\cal L}) = {\cal L}_{\rm ab}$$
If $\cal L$ is an abelian Hom-Leibniz $n$-algebra, then  ${\ele}$ is endowed with a trivial Hom-co-representation structure from ${\cal D}_{n-1}\left( {\cal L} \right)$, then
$${_n}HL_{1}^{\alpha}\left( {\cal L} \right) =HL_{1}^{\alpha}({\cal D}_{n-1}\left( {\cal L} \right), {\cal L})=\frac {{\cal L} \otimes {\cal L}^{\otimes\left( n-1\right)}} {\alpha({\cal L})\otimes[{\cal L}^{\otimes \left( n-1\right)},{\cal L}^{\otimes \left( n-1\right)}]}$$

When ${\cal L}$ is a Hom-Leibniz $2$-algebra, that is, a Hom-Leibniz algebra, then we have that
 $${_2}CL_{\ast}^{\alpha}({\cal L})=CL_{\ast}^{\alpha}({\cal L}, {\cal L}) \cong CL_{\ast+1}^{\alpha}(\cal L)$$
(see the proof of Proposition 3.4 in \cite{CIP1}). Thus ${_2}HL_k^{\alpha}({\cal L}) \cong HL_{k+1}^{\alpha}({\cal L})$, for all $k \geq 1$. In particular, ${_2}HL_0^{\alpha}({\cal L}) \cong HL_{1}^{\alpha}({\cal L}) \cong {\cal L}_{\rm ab}$.
When $\alpha=\id$ the corresponding results for Leibniz $n$-algebras in \cite{Ca1, Ca} are recovered.


 \section{Universal central extensions} \label{uce}
  \begin{De} \label{n-alfacentral}
 A short exact sequence of Hom-Leibniz $n$-algebras $(K) : 0 \to ({\eme}, \widetilde \alpha_{\eme}) \stackrel{i} \to ({\ka},\widetilde \alpha_{\ka}) \stackrel{\pi} \to ({\ele}, \widetilde \alpha_{\ele}) \to 0$ is said to be central if $[{\eme}, {\ka}, \stackrel{n-1} \dots, {\ka}]=0$.
 Equivalently,  ${\eme} \subseteq Z({\ka})$.

We say that $(K)$ is $\alpha$-central if $[\alpha_{\eme}({\eme}), \stackrel{n-1}\dots, \alpha_{\eme}({\eme}), {\ka}]=0$.
\end{De}

\begin{Rem} \label{central}
Let us observe that  the notion of central extension in case $\widetilde{\alpha}_{\cal K}= (\id_{\cal K})$  coincides with the notion of central extension of Leibniz $n$-algebras given in \cite{Ca1}. Nevertheless, the notion of $\alpha$-central extension in case $\widetilde{\alpha}_{\cal K}= (\id_{\cal K})$ gives rise to a new notion of central extension of Leibniz $n$-algebras. In particular, this kind of central extensions are abelian extensions of Leibniz $n$-algebras \cite{CLP}.

In case $n=2$, we recover the notions of central and $\alpha$-central extension of a Hom-Leibniz algebra introduced in \cite{CIP1}.

Obviously every central extension is an $\alpha$-central extension, but the converse doesn't hold as the following counterexample shows:

Let $({\cal L}, \widetilde{\alpha}_{\cal L})$ be the Hom-Leibniz 3-algebra  where ${\cal L}$ is the two-dimensional vector space with basis $\{a_1, a_2\}$, the bracket operation is given by $[a_i,a_i,a_i]=a_i, i=1, 2$ and zero elsewhere, and endomorphism $\alpha_{\cal L} = 0$.

On the other hand, let $({\ka}, \widetilde{\alpha}_{\ka})$ be the Hom-Leibniz 3-algebra  where ${\ka}$ is the three-dimensional vector space with basis $\{b_1, b_2, b_3\}$, the bracket operation is given by $[b_i,b_i,b_i]=b_i, i=1, 2, 3$ and zero elsewhere, and endomorphism $\alpha_{\ka} = 0$.

The surjective homomorphism $\pi : ({\ka}, \widetilde{\alpha}_{\ka}) \twoheadrightarrow ({\ele}, \widetilde{\alpha}_{\ele})$  given by $\pi(b_1)=0,$ $\pi(b_2)=a_1, \pi(b_3)=a_2$, is an $\alpha$-central extension, since $\Ker(\pi) = \langle \{b_1\} \rangle$ and $[\alpha_{\ka}(\Ker(\pi)), \alpha_{\ka}(\Ker(\pi)), {\ka}] = 0$, but is not a central extension since $Z({\ka}) = 0$.
\end{Rem}

\begin{De} \label{Def universal}
A central extension $(K) : 0 \to ({\cal M},\widetilde{\alpha}_{\eme}) \stackrel{i} \to ({\cal K},\widetilde{\alpha}_{\ka}) \stackrel{\pi} \to ({\cal L},\widetilde{\alpha}_{\cal L}) \to $$0$ is said to be universal if for every central extension $(K') : 0 \to ({\cal M}',\widetilde {\alpha}_{{\eme}'}) \stackrel{i'} \to ({\cal K}', \widetilde{\alpha}_{{\ka}'}) \stackrel{\pi'} \to ({\cal L},\widetilde{\alpha}_{\cal L}) \to $$0$ there exists a unique homomorphism of Hom-Leibniz $n$-algebras  $h : ({\ka}, \widetilde{\alpha}_{\ka}) \to ({\ka}',\widetilde{\alpha}_{\ka'})$ such that $\pi' \circ h = \pi$.

The central extension  $(K) : 0 \to ({\cal M},\widetilde{\alpha}_{\cal M}) \stackrel{i} \to ({\cal K},\widetilde{\alpha}_{\cal K}) \stackrel{\pi} \to ({\cal L},\widetilde{\alpha}_{\ele}) \to 0$ is said to be universal $\alpha$-central extension if for every $\alpha$-central extension $(K') : 0 \to ({\cal M'},\widetilde{\alpha}_{{\cal M}'}) \stackrel{i'} \to ({\cal K}',\widetilde{\alpha}_{{\cal K}'}) \stackrel{\pi'} \to ({\cal L},\widetilde{\alpha}_{\cal L}) \to 0$ there exists a unique homomorphism of Hom-Leibniz $n$-algebras  $h : (\cal K,\widetilde \alpha_K) \to (\cal K',\widetilde \alpha_{K'})$ such that $\pi' \circ h = \pi$.
\end{De}

\begin{Rem} \label{rem}
Obviously, every universal $\alpha$-central extension is a universal central extension.
Note that in the case $\widetilde{\alpha}_{\ka} = (\id_{\ka})$  both notions coincide.
In case $n=2$ we recover the corresponding notions of universal ($\alpha$-)central extension  of Hom-Leibniz algebras given respectively in \cite[Definition 4.3]{CIP1}.
\end{Rem}

\begin{De}
A Hom-Leibniz $n$-algebra $({\cal L},\widetilde{\alpha}_{\cal L})$ is said to be perfect if ${\cal L} = [{\cal L}, \dots, {\cal L}]$.
\end{De}

\begin{Le} \label{lema 1}
Let $\pi : ({\cal K},\widetilde{\alpha}_{\ka}) \twoheadrightarrow ({\cal L},\widetilde{\alpha}_{\cal L})$ be a surjective homomorphism of Hom-Leibniz $n$-algebras. If $({\cal K},\widetilde{\alpha}_{\ka})$ is a perfect  Hom-Leibniz $n$-algebra, then $({\cal L},\widetilde{\alpha}_{\ele})$  is a perfect Hom-Leibniz $n$-algebra as well.
\end{Le}

\begin{Le} \label{lema 3}
If $0 \to \left( {\cal M},\widetilde{\alpha }_{\cal M}\right) \stackrel{i}\to \left( {\cal K},\widetilde{\alpha }_{\cal K}\right) \stackrel{\pi}\to \left( {\cal L},\widetilde{\alpha }_{\cal L}\right) \to 0$ is a universal central extension, then $\left( {\cal K},\widetilde{\alpha }_{\cal K}\right)$ and $\left( {\cal L},\widetilde{\alpha }_{\cal L}\right)$ are perfect Hom-Leibniz $n$-algebras.
\end{Le}
{\it  Proof.} Assume that $({\cal K},\widetilde{\alpha }_{\cal K})$ is not a perfect Hom-Leibniz $n$-algebra, then $\left[ {\cal K}, \dots,{\cal K}\right]$ $\nsubseteq {\cal K}$, thus  $\left( {\cal K}/\left[ {\cal K}, \dots, {\cal K}\right], \widetilde{\alpha} _{\mid }\right)$, where $(\widetilde{\alpha} _{\mid})$ is the induced natural homomorphism, is an abelian Hom-Leibniz $n$-algebra (see Example \ref{ejemplo 1} {\it c)}).

Consider the central extension $0\to ({\cal K}_{\rm ab}, \widetilde{\alpha}_{\mid})\to \left( {\cal K}_{\rm ab} \times {\cal L},\widetilde{\alpha}_{\mid} \times \widetilde{\alpha}_{\cal L}\right) \overset{pr}\to \left( {\cal L},\widetilde{\alpha }_{\cal L}\right)$ $\to 0$, then the homomorphisms of Hom-Leibniz $n$-algebras $\varphi ,\psi :({\cal K},\widetilde{\alpha}_{\ka})\to \left( {\cal K}_{\rm ab} \times {\cal L}, \right.$ $\left.\widetilde{\alpha}_{\mid} \times \widetilde{\alpha}_{\cal L}\right)$ given by
$\varphi \left( k\right) =\left( \overline{k},\pi \left( k\right) \right)$ and $\psi \left( k\right) =\left( 0,\pi \left( k\right) \right)$ ($\overline{k}$ denotes the coset $k+[{\cal K}, \dots, {\cal K}]$) are two distinct homomorphisms of Hom-Leibniz $n$-algebras such that $pr\circ \varphi =\pi =pr\circ \psi$, which contradicts the universality of the given extension.

Lemma \ref{lema 1} completes the proof. \rdg

\begin{Le}\label{lema 2}
 Let $0 \to ({\cal M},\widetilde{\alpha}_{\cal M}) \stackrel{i} \to ({\cal K},\widetilde{\alpha}_{\ka}) \stackrel{\pi} \to ({\cal L},\widetilde{ \alpha}_{\cal L}) \to 0$ be an $\alpha$-central extension and $({\ka},\widetilde{\alpha}_{\ka})$ is a  perfect  Hom-Leibniz $n$-algebra. If there exists a homomorphism of Hom-Leibniz $n$-algebras $f : ({\cal K},\widetilde{\alpha}_{\cal K}) \to ({\cal A},\widetilde{\alpha}_{\as})$ such that $\tau \circ f = \pi$, where $ 0 \to ({\cal N},\widetilde{\alpha}_{\cal N}) \stackrel{j} \to ({\cal A},\widetilde{\alpha}_{\cal A}) \stackrel{\tau} \to ({\cal L},\widetilde{\alpha}_{\cal L}) \to 0$ is a  central  extension, then $f$ is unique.
\end{Le}
{\it Proof.} Assume that there are two homomorphisms $f_{1},f_{2}:({\ka},\widetilde{\alpha}_{\ka})\to ({\as},\widetilde{\alpha}_{\as})$ such that $\tau \circ f_{1}= \pi = \tau \circ f_{2}$, then $f_{1}-f_{2}\in \Ker(\tau) = {\ene}$, i.e. $f_{1}(k)=f_{2}(k)+n_{k},$ $n_{k}\in {\ene}$.

 Since $({\ka},\widetilde{\alpha}_{\ka})$ is a perfect Hom-Leibniz $n$-algebra, it is enough to show that  $f_{1}$ and $f_{2}$ coincide on $\left[ {\ka},\dots ,{\ka} \right]$. Indeed
\[\begin{array}{l}
f_{1}\left[  k_{1}, \dots, k_{n}\right] = \left[  f_{2}\left( k_{1}\right)  + n_{k_{1}}, \dots, f_{2}\left( k_{n} \right)+ n_{k_{n}} \right] = \left[ f_{2}\left( k_{1}\right), \dots, f_{2}\left( k_{n}\right) \right] + A =\\
 f_{2}\left[  k_{1}, \dots, k_{n} \right],
\end{array}\]
\noindent since a typical summand  in $A$ is of the form $\left[ n_{k_{1}}, \dots,n_{{k}_j}, f_{2}\left( k_{j+1}\right),\dots, f_{2}\left( k_{n}\right) \right]$ which vanishes because ${\ene} \subseteq Z({\as})$. \rdg

\bigskip

The category ${\sf {_n}HomLeib}$ is a semi-abelian category that doesn't satisfy the so called in \cite{CVdL} {\bf UCE} condition, namely the composition of central extensions is a central extension, as the following example shows:

\begin{Ex} \label{uce condition}
Let $\left( {\cal L},\widetilde{\alpha }_{\cal L}\right) $ be the two-dimensional Hom-Leibniz $3$-algebra with basis $\left\{ b_{1},b_{2}\right\}$, bracket given by $\left[ b_{2},b_{1},b_{1}\right] =b_{2}, \left[ b_{2},b_{2},b_{2}\right] =b_{1}$ and zero elsewhere, and endomorphism $\widetilde{\alpha }_{\cal L}=(0).$

 Let $\left( {\cal K},\widetilde{\alpha }_{\cal K}\right)$ be the three-dimensional Hom-Leibniz $3$-algebra with basis $\left\{ a_{1},a_{2},\right.$ $\left.a_{3}\right\}$, bracket given by $\left[ a_{2},a_{2},a_{2}\right] =a_{1}, \left[ a_{3},a_{2},a_{2}\right]
=a_{3},\left[ a_{3},a_{3},a_{3}\right] =a_{2}$ and zero elsewhere, and endomorphism $\widetilde{\alpha }_{\cal K}=(0).$

Obviously $\left( {\cal K},\widetilde{\alpha }_{\cal K}\right)$ is a perfect Hom-Leibniz $3$-algebra and $Z({\cal K})=\langle \left\{ a_{1}\right\} \rangle.$
The linear map $\pi :\left( {\cal K},\widetilde{0}\right) \to \left( {\ele},\widetilde{0}\right)$ given by $\pi \left( a_{1}\right) =0, \pi \left( a_{2}\right) =b_{1}, \pi \left( a_{3}\right) =b_{2}$, is a central central extension since $\pi$ is a surjective  homomorphism of Hom-Leibniz $3$-algebras and  $\Ker\left( \pi \right) = \langle \left\{ a_{1}\right\} \rangle \subseteq Z({\cal K})$.

Now consider the four-dimensional Hom-Leibniz $3$-algebra $\left( {\cal F},\widetilde{\alpha }_{\cal F}\right)$ with basis $\left\{
e_{1},e_{2},e_{3},e_{4}\right\}$, bracket given by $\left[ e_{3},e_{2},e_{2}\right] =e_{1\text{ }}, \left[ e_{3},e_{3},e_{3} \right] =e_{2}, \left[ e_{4},e_{3},e_{3}\right] =e_{4\text{ }}, \left[ e_{4},e_{4},e_{4}\right] =e_{3}$ and zero elsewhere, and endomorphism $\widetilde{\alpha }_{\cal F}=(0).$

The linear map $\rho \left( e_{1}\right) =0, \rho \left( e_{2}\right)=a_{1}, \rho \left( e_{3}\right) =a_{2}, \rho \left( e_{4}\right) =a_{3}$ is a central extension since $\rho$ is a surjective  homomorphism of Hom-Leibniz $3$-algebras and $\Ker\left( \rho \right) = \langle \left\{ e_{1}\right\} \rangle = Z({\cal F})$.

The composition $\pi \circ \rho :\left( {\cal F}, \widetilde{0}\right) \longrightarrow \left({\cal L}, \widetilde{0}\right)$ is given by $\pi \circ \rho \left( e_{1}\right) =\pi \left( 0\right) =0, \pi \circ \rho \left( e_{2}\right) =\pi \left( a_{1}\right)
=0, \pi \circ \rho \left( e_{3}\right) =\pi \left( a_{2}\right) =b_{1}, \pi \circ \rho \left( e_{4}\right) =\pi \left( a_{3}\right) =b_{2}$. Consequently, $\pi \circ \rho $ is a surjective homomorphism, but is not a central extension, since $\Ker\left( \pi \circ
\rho \right) =\langle \{e_1, e_2\}\rangle  \nsubseteq Z({\cal F})$. However, $\pi \circ \rho :\left( {\cal F}, \widetilde{0}\right) \longrightarrow \left({\cal L}, \widetilde{0}\right)$ is an $\alpha$-central extension.
\end{Ex}

\begin{Le} \label{lema 4}
Let $0 \to ({\cal M} ,\widetilde{\alpha}_{\cal M}) \stackrel{i} \to ({\cal K}, \widetilde{\alpha}_{\cal K}) \stackrel{\pi} \to ({\cal L}, \widetilde{\alpha}_{\cal L}) \to 0$  and $0 \to ({\cal N}, \widetilde{\alpha}_{\cal N}) \stackrel{j} \to ({\cal F},\widetilde{\alpha}_{\cal F}) \stackrel{\rho} \to ({\cal K},\widetilde{\alpha}_{\cal K}) \to 0$ be central extensions with  $({\cal K},\widetilde{\alpha}_{\cal K})$ a perfect Hom-Leibniz $n$-algebra. Then the composition extension $0 \to ({\cal P}, \widetilde{\alpha}_{\cal P}) = \Ker(\pi \circ \rho)   \to ({\cal F}, \widetilde{\alpha}_{\cal F}) \stackrel{\pi \circ \rho} \to ({\cal L}, \widetilde{\alpha}_{\cal L}) \to $$0$ is an $\alpha$-central extension.

Moreover, if $0 \to ({\cal M} ,\widetilde{\alpha}_{\cal M}) \stackrel{i} \to ({\cal K}, \widetilde{\alpha}_{\cal K}) \stackrel{\pi} \to ({\cal L}, \widetilde{\alpha}_{\cal L}) \to 0$ is a universal $\alpha$-central extension, then $0 \to ({\cal N}, \widetilde{\alpha}_{\cal N}) \stackrel{j} \to ({\cal F},\widetilde{\alpha}_{\cal F}) \stackrel{\rho} \to ({\cal K},\widetilde{\alpha}_{\cal K}) \to 0$ is split.
\end{Le}
{\it Proof.} Since $\left( {\cal K},\widetilde{\alpha }_{\ka}\right)$ is a perfect Hom-Leibniz $n$-algebra, then  every element $f\in {\cal F}$ can be written as $\underset{k}{\sum }\lambda _{k}\left[ f_{k1},...,f_{kn}\right]+n,$ $n\in {\cal N},$ $f_{k1},...,f_{kn}\in {\cal F}$. So, for any element in $[\alpha_{\cal P}({\cal P}), \stackrel{n-1} \dots, \alpha_{\cal P}({\cal P}), {\cal F}]$ we have
\[
\begin{array}{l}
\left[\alpha_{\cal P} \left(p_{1}\right), \dots, f_{i}, \dots, \alpha_{\cal P} \left( p_{n-1}\right) \right] =\\
\underset{k}{\sum }\lambda _{k}\left( \left[\alpha_{\cal P} \left(p_{1}\right), \dots,  [f_{i_{k1}}, \dots, f_{i_{kn}}], \dots, \alpha_{\cal P} \left( p_{n-1}\right) \right] \right. +\\
\left. \left[\alpha_{\cal P} \left(p_{1}\right), \dots, n, \dots, \alpha_{\cal P} \left( p_{n-1}\right) \right] \right)
\end{array}
\]
which vanishes by application of the fundamental identity and bearing in mind that $[{\cal P},{\cal F}, \dots, {\cal F}] \in \Ker(\rho) = {\ene}$  and ${\cal N} \subseteq Z({\cal F})$.

For the second statement, if $0 \to ({\cal M} ,\widetilde{\alpha}_{\cal M}) \stackrel{i} \to ({\cal K}, \widetilde{\alpha}_{\cal K}) \stackrel{\pi} \to ({\cal L}, \widetilde{\alpha}_{\cal L}) \to 0$ is a universal $\alpha$-central extension, then by the first statement, $0 \to ({\cal P}, \widetilde{\alpha}_{\cal P}) = \Ker(\pi \circ \rho)   \to ({\cal F}, \widetilde{\alpha}_{\cal F}) \stackrel{\pi \circ \rho} \to ({\cal L}, \widetilde{\alpha}_{\cal L}) \to $$0$ is an $\alpha$-central extension, then there exists a unique homomorphism of Hom-Leibniz algebras $\sigma : ({\cal K}, \widetilde{\alpha}_{\cal K}) \to ({\cal F}, \widetilde{\alpha}_{\cal F})$ such that $\pi \circ \rho \circ\sigma = \pi$. On the other hand, $\pi \circ \rho \circ \sigma = \pi = \pi \circ \id$ and $({\cal K}, \widetilde{\alpha}_{\cal K})$ is perfect, then Lemma \ref{lema 2} implies that $\rho \circ \sigma = \id$. \rdg

\begin{Th}\label{teorema}\
\begin{enumerate}
\item[a)]  If a central extension $0 \to ({\cal M} ,\widetilde{\alpha}_{\cal M}) \stackrel{i} \to ({\cal K}, \widetilde{\alpha}_{\cal K}) \stackrel{\pi} \to ({\cal L}, \widetilde{\alpha}_{\cal L}) \to 0$ is a universal $\alpha$-central extension, then
$({\cal K}, \widetilde{\alpha}_{\cal K})$ is a perfect Hom-Leibniz $n$-algebra and every central
extension of $({\cal K}, \widetilde{\alpha}_{\cal K})$ is split.

    \item[b)] Let  $0 \to ({\eme}, \widetilde{\alpha}_{\eme}) \stackrel{i} \to ({\ka},\widetilde{\alpha}_{\ka}) \stackrel{\pi} \to ({\ele}, \widetilde{\alpha}_{\ele}) \to 0$ be a central extension.

     If $({\ka},\widetilde{\alpha}_{\ka})$ is a perfect Hom-Leibniz $n$-algebra and every central extension of $({\ka},\widetilde{\alpha}_{\ka})$ is split, then $0 \to ({\eme}, \widetilde{\alpha}_{\eme}) \stackrel{i} \to ({\ka},\widetilde{\alpha}_{\ka}) \stackrel{\pi} \to ({\ele}, \widetilde{\alpha}_{\ele}) \to 0$ is a universal central extension.

\item[c)] A Hom-Leibniz $n$-algebra  $({\ele}, \widetilde{\alpha}_{\ele})$ admits a universal central extension if and only if  $({\ele}, \widetilde{\alpha}_{\ele})$ is perfect. Furthermore, the kernel of the universal central extension is canonically isomorphic to ${_n}HL_{1}^{\alpha}\left( {\cal L} \right)$.


\item[d)] If $0 \to ({\eme}, \widetilde{\alpha}_{\eme}) \stackrel{i} \to ({\ka},\widetilde{\alpha}_{\ka}) \stackrel{\pi} \to ({\ele}, \widetilde{\alpha}_{\ele}) \to 0$ is a universal $\alpha$-central extension, then ${_n}HL_{0}^{\alpha}\left( {{\ka}} \right) = {_n}HL_{1}^{\alpha}\left( {{\ka}} \right) = 0$.

    \item[e)] If ${_n}HL_{0}^{\alpha}\left( {\ka} \right) = {_n}HL_{1}^{\alpha}\left( {\ka} \right) = 0$, then any central extension $0 \to ({\eme}, \widetilde{\alpha}_{\eme}) \stackrel{i} \to ({\ka}, \widetilde{\alpha}_{\ka}) \stackrel{\pi} \to ({\ele}, \widetilde{\alpha}_{\ele}) \to 0$ is a universal central extension.
\end{enumerate}
\end{Th}
{\it Proof.}\
\noindent {\it a)} If $0 \to ({\eme},\widetilde{\alpha}_{\eme}) \stackrel{i} \to ({\ka}, \widetilde{\alpha}_{\ka}) \stackrel{\pi} \to ({\ele}, \widetilde{\alpha}_{\ele}) \to 0$ is a universal $\alpha$-central extension, then is a universal central extension by Remark \ref{rem}, so $({\ka}, \widetilde{\alpha}_{\ka})$ is a perfect Hom-Leibniz $n$-algebra by Lemma \ref{lema 3} and every central extension of $({\ka}, \widetilde{\alpha}_{\ka})$ is split by Lemma \ref{lema 4}.

\bigskip

\noindent {\it b)} Let us consider  any central extension $0 \to ({\ene}, \widetilde{\alpha}_{\ene}) \stackrel{j} \to ({\as},\widetilde{\alpha}_{\as}) \stackrel{\tau} \to ({\ele}, \widetilde{\alpha}_{\ele}) \to 0$. Construct the pull-back extension  $0 \to ({\ene}, \widetilde{\alpha}_{\ene}) \stackrel{\chi} \to ({\qu},\widetilde{\alpha}_{\qu}) \stackrel{\overline{\tau}} \to ({\ka}, \widetilde{\alpha}_{\ka}) \to 0$, where ${\qu} = {\as} \times_{\ele}  {\ka} =\{(a,k) \in {\as} \times {\ka} \mid \tau(a)=\pi(k) \}$ and $\alpha_{\qu}(a,k)=(\alpha_{\as}(a),\alpha_{\ka}(k))$, which is central, consequently is split, that is, there exists a homomorphism $\sigma : ({\ka},\widetilde{\alpha}_{\ka}) \to ({\qu},\widetilde{\alpha}_{\qu})$ such that $\overline{\tau} \circ \sigma = \id$.

Then $\overline{\pi} \circ \sigma$, where $\overline{\pi} :  ({\qu},\widetilde{\alpha}_{\qu})\to ({\as},\widetilde{\alpha}_{\as})$ is induced by the pull-back construction, satisfies $\tau \circ \overline{ \pi} \circ \sigma = \pi$. Lemma \ref{lema 2} ends the proof.

\bigskip

\noindent {\it c)} For a Hom-Leibniz $n$-algebra $\left( {\ele}, \widetilde{\alpha} _{\ele} \right)$, consider the chain homology complex $_{n}C_{\ast }^{\alpha }\left( {\ele},\mathbb{K}\right)$, where $\mathbb{K}$ is endowed with a trivial Hom-co-representation structure.
$$_{n}C_{\ast }^{\alpha }\left( {\ele},\mathbb{K}\right):  \dots \to {\ele}^{\otimes k(n-1)+1} \stackrel{\delta_k}\to {\ele}^{\otimes (k-1)(n-1)+1} \stackrel{\delta_{k-1}} \to \dots \to {\ele}^{\otimes 2n-1} \stackrel{\delta_2}\to {\ele}^{\otimes n} \stackrel{\delta_1}\to {\ele}$$
 The low differentials are given by
 \[\begin{array}{l}
 \delta_1(x_1 \otimes \dots \otimes x_n) = [x_1, \dots, x_n]\\
\delta_2(x_1 \otimes \dots \otimes x_n \otimes y_1 \otimes \dots \otimes y_{n-1}) = [x_1, \dots, x_n]\otimes \alpha_{\ele}(y_1) \otimes \dots \otimes \alpha_{\ele}(y_{n-1}) -\\
 \displaystyle \sum_{i=1}^n \alpha_{\ele}(x_1) \otimes \dots \otimes [x_i,y_1, \dots,y_{n-1}] \otimes \dots \otimes \alpha_{\ele}(x_{n})
\end{array} \]

As a $\mathbb{K}$-vector space, let $I_{\ele}$ be the subspace of ${\ele}^{\otimes 2n-1}$ spanned
by the elements of the form
\[
\begin{array}{l} \left[ x_{1}, \dots,x_{n}\right] \otimes \alpha _{L}\left( y_{1}\right) \otimes \dots \otimes \alpha _{L}\left( y_{n-1}\right) -\\
 \displaystyle \sum_{i=1}^n \alpha_{\ele}(x_1) \otimes \dots \otimes [x_i,y_1, \dots,y_{n-1}] \otimes \dots \otimes \alpha_{\ele}(x_{n}) \end{array} \]
that is $I_{\ele}=\im \left( \delta _{2}:{\ele}^{\otimes 2n-1}\to {\ele}^{\otimes n}\right)$. Let $\frak{uce}(\ele)$ be the quotient
vector space $\frac{{\ele}^{\otimes n}}{I_{\ele}}$. Every coset $\left( x_{1}\otimes \dots \otimes x_{n}\right) +I_{\ele}$ is denoted by $\left\{ x_{1}, \dots,x_{n}\right\}$.

By construction, the following identity holds
\begin{equation} \label{identidad}
\begin{array}{l} \left\{ \left[ x_{1}, \dots,x_{n}\right], \alpha _{L}\left( y_{1}\right), \dots, \alpha _{L}\left( y_{n-1}\right) \right\} =\\
\displaystyle \sum_{i=1}^n \left\{ \alpha_{\ele}(x_1), \dots, [x_i,y_1, \dots,y_{n-1}], \dots, \alpha_{\ele}(x_{n}) \right\} \end{array}
\end{equation}

Since $\delta _{1}$ vanishes on $I_{\ele}$, it induces a linear map $u_{\ele}: \frak{uce}({\ele}) \to {\ele}$, given by $u_{\ele}(\left\{ x_{1}, \dots,x_{n}\right\}) =\left[ x_{1}, \dots,x_{n}\right]$, and $(\widetilde{\alpha}_{\ele})$ induces $(\widetilde{\alpha}_{\frak{uce}(\ele)})$, where $\alpha_{\frak{uce}(\ele)}(\{x_1, \dots, x_n\})$ $= \{\alpha_{\ele}(x_1), \dots, \alpha_{\ele}(x_n)\}$.

The bracket operation
 $$\left[ \left\{ x_{1,1},...,x_{n,1}\right\}, \dots,\left\{x_{1,n}, \dots,x_{n,n}\right\} \right] =\left\{ \left[ x_{1,1}, \dots,x_{n,1} \right], \dots,\left[ x_{1,n}, \dots ,x_{n,n}\right] \right\}$$
endows  $\left( \frak{uce}({\ele}), \widetilde{\alpha}_{\frak{uce}(\ele)} \right)$  with a structure of  Hom-Leibniz $n$-algebra
and $u_{\ele}:\left( \frak{uce}({\ele}), \widetilde{\alpha}_{\frak{uce}(\ele)} \right) \to  \left( {\ele},\widetilde{\alpha }_{\ele}\right)$ becomes an epimorphism of Hom-Leibniz $n$-algebras when $\left( {\ele}, \widetilde{\alpha} _{\ele} \right)$ is perfect because $\im(u_{\ele})=[{\ele}, \dots, {\ele}]$.

From the construction immediately follows that $\Ker\left( u_{\ele}\right)={_n}HL_{1}^{\alpha}({\ele})$, so we have the central extension
$$0\to \left( _{n}HL_{1}^{\alpha }\left( {\ele} \right), \widetilde{\alpha}_{\frak{uce}(\ele)\mid}\right) \to \left( \frak{uce}({\ele}), \widetilde{\alpha}_{\frak{uce}(\ele)} \right) \stackrel{u_{\ele}}\to  \left( {\ele},\widetilde{\alpha }_{\ele}\right) \to 0$$
which is universal, because for any central extension $0\to \left( {\eme}, \widetilde{\alpha }_{\eme}\right) \to
\left( {\ka},\widetilde{\alpha }_{\ka}\right) \overset{\pi} \to \left( {\ele}, \widetilde{\alpha }_{\ele}\right) \to 0$ there exists the homomorphism of Hom-Leibniz $n$-algebras $\beta :\left( \frak{uce}({\ele}), \widetilde{\alpha}_{\frak{uce}(\ele)} \right)  \to
\left( {\ka},\widetilde{\alpha }_{\ka}\right)$ given by $\beta \left( \left\{ x_{1}, \dots, x_{n}\right\} \right) =\left[ k_{1},...,k_{n}\right] ,\pi \left( k_{i}\right) =x_{i}$, such that $\pi \circ \beta =u_{\ele}$.

A direct checking shows that $\left( \frak{uce}({\ele}), \widetilde{\alpha}_{\frak{uce}(\ele)} \right)$ is perfect, then Lemma \ref{lema 2} guarantees the uniqueness of $\beta$.

\bigskip

\noindent {\it d)} If $0\to \left( {\eme}, \widetilde{\alpha }_{\eme}\right) \to \left( {\ka},\widetilde{\alpha }_{\ka}\right) \overset{\pi} \to \left( {\ele}, \widetilde{\alpha }_{\ele}\right) \to 0$ is a universal $\alpha$-central extension, then $({\ka} ,\widetilde{\alpha}_{\ka})$ is perfect by Remark \ref{rem} and Lemma \ref{lema 3}, so  ${_n}HL_{0}^{\alpha}\left( {\ka} \right) = 0$. By Lemma \ref{lema 4} and statement {\it c)}, the universal central extension corresponding to $({\ka}, \widetilde{\alpha}_{\ka})$ is split, so $ {_n}HL_{1}^{\alpha}\left( {\ka} \right) = 0$.

\bigskip

\noindent {\it e)} \ ${_n}HL_{0}^{\alpha}\left( {\cal K} \right) = 0$ implies that $({\cal K} ,\widetilde{\alpha}_{\cal K})$ is a
perfect Hom-Leibniz $n$-algebra.

${_n}HL_{1}^{\alpha}\left( {\ka} \right) = 0$ implies that $(\frak{uce}(\ka),\widetilde{\alpha}_{\frak{uce}({\ka})}) \stackrel{\sim} \to ({\ka} ,\widetilde{\alpha}_{\ka})$. Statement {\it b)} ends the proof.
 \rdg

\begin{Rem}
 When $n=2$, the above results recover the corresponding ones for Hom-Leibniz algebras in \cite{CIP1}.
\end{Rem}


\section{Non-abelian tensor product} \label{non-abelian tensor}

Let $({\eme}_i,\widetilde{\alpha}_{{\eme}_i}), 1 \leq i \leq n$,  be $n$-sided Hom-ideals of a Hom-Leibniz $n$-algebra $({\ele}, \widetilde{\alpha}_{{\ele}})$. We denote by ${\eme}_1 \ast \dots \ast {\eme}_n$  the vector space spanned by all the symbols $m_{\sigma(1)} \ast \dots \ast m_{\sigma(n)}$, where $m_i \in {\eme}_i, i \in \{1, 2, \dots, n\}, \sigma \in S_n$.

We claim that $({\eme}_1 \ast \dots \ast {\eme}_n, \widetilde{\alpha}_{{\eme}_1 \ast \dots \ast {\eme}_n})$ is a Hom-vector space, where $(\widetilde{\alpha}_{{\eme}_1 \ast \dots \ast {\eme}_n})$ is  induced by $\alpha_{{\eme}_i}, 1 \leq i \leq n$, i.e.
\begin{align*}
\alpha_{{\eme}_1 \ast \dots \ast {\eme}_n}\left( m_{\sigma(1)} \ast \dots \ast m_{\sigma(n)} \right)  =\alpha_{{\eme}_{\sigma(1)}}\left(m_{\sigma(1)} \right)  \ast \dots \ast \alpha_{{\eme}_{\sigma(n)}} \left( m_{\sigma(n)} \right).
\end{align*}

We denote by ${\cal DL}_n({\eme}_1, \dots, {\eme}_n)$ the vector subspace  spanned by the elements of the form:

\begin{enumerate}
\item[{\it a)}] $\lambda \left( m_{\sigma (1)}\ast \dots \ast m_{\sigma(n)}\right) =\left( \lambda m_{\sigma (1)}\right) \ast m_{\sigma (2)}\ast \dots \ast  m_{\sigma (n)}= \dots = m_{\sigma (1)}\ast \dots \ast  \left(\lambda  m_{\sigma (n)}\right)$.

\item[{\it b)}] $m_{\sigma (1)}\ast \dots \ast \left( m'_{\sigma (i)} + m''_{\sigma (i)} \right) \ast \dots \ast m_{\sigma(n)} =
m_{\sigma (1)}\ast \dots \ast m'_{\sigma (i)}  \ast \dots \ast m_{\sigma(n)} + m_{\sigma (1)}\ast \dots \ast  m''_{\sigma (i)} \ast \dots \ast m_{\sigma(n)}$, for any  $i \in \{1, 2, \dots, n\}$.

\item[{\it c)}] $\left[ m_{\tau (1)}, \dots, m_{\tau(n)}\right] \ast \alpha _{{\eme}_{\tau (n+1)}}\left( m_{\tau (n+1)}\right) \ast \dots \ast \alpha_{{\eme}_{\tau(2n-1)}}\left( m_{\tau(2n-1)}\right) -$

$\displaystyle \sum_{i=1}^n \alpha _{{\eme}_{\tau(1)}}(m_{\tau (1)}) \ast \dots \ast \left[  m_{\tau (i)},  m_{\tau (n+1)}, \dots,  m_{\tau (2n-1)}  \right] \ast \dots \ast \alpha _{{\eme}_{\tau(n)}}(m_{\tau (n)})$.

\item[{\it d)}] $[m_{\sigma(1)}, \dots, m_{\sigma(n)}] \ast \alpha_{{\eme}_{n+1}}(m_{n+1}) \ast \dots \ast \alpha_{{\eme}_{2n-1}}(m_{2n-1}) - (-1)^{\epsilon(\sigma)} [m_1, \dots, m_n] \ast \alpha_{{\eme}_{n+1}}(m_{n+1}) \ast \dots \ast \alpha_{{\eme}_{2n-1}}(m_{2n-1})$.
\end{enumerate}
for all $\lambda \in \K$, $m_{i} \in {\eme}_i, 1 \leq i \leq n$, $\sigma \in S_n, \tau \in S_{2n-1}$.

Moreover, it can be readily checked  that  $\alpha_{{\eme}_1 \ast \dots \ast {\eme}_n}({\cal DL}_n({\eme}_1, \dots, {\eme}_n)) \subseteq {\cal DL}_n({\eme}_1, \dots, {\eme}_n)$, hence we can construct the quotient Hom-vector space
$$\left({\eme}_1 \ast \dots \ast {\eme}_n / {\cal DL}_n({\eme}_1, \dots, {\eme}_n), \overline{\alpha}_{{\eme}_1 \ast \dots \ast {\eme}_n}\right)$$
which is endowed with a structure of Hom-Leibniz $n$-algebra with respect to the bracket
\begin{equation} \label{bracket}
\begin{array}{l}
[m_{11} \ast \dots \ast m_{n1}, m_{12} \ast \dots \ast m_{n2}, \dots, m_{1n} \ast \dots \ast m_{nn}]:= \\

[m_{11}, \dots, m_{n1}] \ast [m_{12}, \dots, m_{n2}] \ast \dots \ast  [m_{1n}, \dots, m_{nn}]
\end{array}
\end{equation}
where we abbreviate a coset $\overline{m_{1i} \ast \dots \ast m_{ni}}$ by   $m_{1i} \ast \dots \ast m_{ni}$ and the endomorphism $\overline{\alpha}_{{\eme}_1 \ast \dots \ast {\eme}_n}$ by $\alpha_{{\eme}_1 \ast \dots \ast {\eme}_n}$.

\begin{De} \label{non abelian tensor}
The above Hom-Leibniz  $n$-algebra structure on $$\left({\eme}_1 \ast \dots \ast {\eme}_n / {\cal DL}_n({\eme}_1, \dots, {\eme}_n), \overline{\alpha}_{{\eme}_1 \ast \dots \ast {\eme}_n}\right)$$ is called the non-abelian tensor product of the $n$-sided Hom-ideals
$({\eme}_i,\widetilde{\alpha}_{{\eme}_i}), 1 \leq i \leq n$, and it will be denoted by $\left({\eme}_1 \ast \dots \ast {\eme}_n , \widetilde{\alpha}_{{\eme}_1 \ast \dots \ast {\eme}_n}\right)$.
\end{De}

\begin{Rem}
If $\widetilde{\alpha}_{\ele}=(\id_{\ele})$, then $\left({\eme}_1 \ast \dots \ast {\eme}_n , \widetilde{\alpha}_{{\eme}_1 \ast \dots \ast {\eme}_n}\right)$ coincides with the non-abelian tensor product of Leibniz $n$-algebras introduced in \cite{Ca2}. In case $n=2$, we recover a particular case of the non-abelian tensor product of Hom-Leibniz algebras given in \cite{CKP1}.
\end{Rem}

For any $n$-sided Hom-ideals $({\eme}_i,\widetilde{\alpha}_{{\eme}_i}), 1 \leq i \leq n$,   of a Hom-Leibniz $n$-algebra $({\ele}, \widetilde{\alpha}_{{\ele}})$, there exists a homomorphism of Hom-Leibniz $n$-algebras
$$\psi : \left({\eme}_1 \ast \dots \ast {\eme}_n , \widetilde{\alpha}_{{\eme}_1 \ast \dots \ast {\eme}_n}\right) \to \left( \displaystyle \bigcap_{i=1}^n {\eme}_i, \widetilde{\alpha}_{\cap} \right)$$
given by $$\psi(m_1 \ast \dots \ast m_n) = [m_1, \dots, m_n]$$
In particular, when ${\eme}_i = {\ele}, 1 \leq i \leq n$, from relation (\ref{bracket}) immediately follows that $\psi : \left({\ele} \ast \dots \ast {\ele}, \widetilde{\alpha}_{{\ele} \ast \dots \ast {\ele}}\right) \twoheadrightarrow \left( [{\ele}, \dots, {\ele}], \widetilde{\alpha}_{\ele \mid} \right)$ is a central extension.

\begin{Th} \label{uce tensor}
If $({\ele}, \widetilde{\alpha}_{{\ele}})$ is a perfect Hom-Leibniz $n$-algebra, then
$\psi : \left({\ele} \ast \dots \ast {\ele}, \right.$ $\left. \widetilde{\alpha}_{{\ele} \ast \dots \ast {\ele}}\right) \twoheadrightarrow \left( {\ele}, \widetilde{\alpha}_{\ele} \right)$
is a universal central extension.
\end{Th}
{\it Proof.}    Let $0\to \left( Ker\left( \chi \right) ,\widetilde{\alpha }_{{\mathcal{C\mid }}}\right) \overset{i}{\to }({\mathcal{C}},\widetilde{\alpha }_{{\mathcal{C}}})\overset{\chi }{\to }\left( {\mathcal{L}},\widetilde{\alpha}_{{\mathcal{L}}}\right) \to 0$ be a central extension of $\left( {\mathcal{L}},\widetilde{\alpha }_{{\mathcal{L}}}\right)$.

Since $\Ker\left( \chi \right) \subseteq Z\left( {\mathcal{C}}\right)$ we get a well-defined homomorphism of Hom-Leibniz $n$-algebras $f:{\mathcal{L}}\ast \dots \ast {\mathcal{L\to C}}$ given on generators by $f\left( l_{1}\ast \dots \ast l_{n}\right) =\left[ c_{l_{1}},...,c_{l_{n}}\right] $, where $c_{l_{i}}$ is an element in $\chi^{-1}\left( l_{i}\right)$,  $i=1,...,n$.

On the other hand, relation (\ref{bracket}) implies that $\left( {\ele} \ast \dots \ast {\ele}, \widetilde{\alpha}_{{\ele} \ast \dots \ast {\ele}} \right)$ is perfect, then the homomorphism $f$ is unique by Remark \ref{central} and  Lemma \ref{lema 2}. \rdg

\begin{Rem} \label{nucleo ext central Leib}
If $({\ele}, \widetilde{\alpha}_{\ele})$ is a perfect Hom-Leibniz $n$-algebra, then $\Ker(\psi) \cong  {_n}HL_1^{\alpha}({\ele})$ by Theorem \ref{teorema} {\it c)}.

Since universal central extensions of perfect Hom-Leibniz $n$-algebras are unique up to isomorphisms, then $\left({\ele} \ast \dots \ast {\ele}, \widetilde{\alpha}_{{\ele} \ast \dots \ast {\ele}}\right) \cong \left( \frak{uce}({\ele}), \widetilde{\alpha}_{\frak{uce}({\ele})} \right)$ by means of the isomorphism $\varphi : \left({\ele} \ast \dots \ast {\ele}, \widetilde{\alpha}_{{\ele} \ast \dots \ast {\ele}}\right) \to \left( \frak{uce}({\ele}), \widetilde{\alpha}_{\frak{uce}({\ele})} \right), \varphi(l_1 \ast \dots \ast l_n) = \{l_1, \dots, l_n\}.$

In case $n=2$, the universal central extension in Theorem \ref{uce tensor} provides the universal central extension of a Hom-Leibniz algebra given in \cite{CKP1}.
\end{Rem}

\begin{Pro}\label{exact-tensor-2}
 If $({\eme},\widetilde{\alpha}_{\eme})$ is an $n$-sided Hom-ideal of a perfect  Hom-Leibniz $n$-algebra $({\ele}, \widetilde{\alpha}_{\ele})$, then there is an exact sequence of vector spaces
 \footnotesize{ \[
\Ker ( \displaystyle \bigoplus_{i=1}^n {\ele} \ast \dots \ast \overbrace{{\eme}}^i \ast \dots \ast {\ele} \stackrel{\psi_{\mid}}\to {\eme}  ) \to {_n}HL_1^{\alpha}({\ele}) \to {_n}HL_1^{\alpha}({\ele}/{\eme}) \to \frac{{\eme}}{\displaystyle \bigoplus_{i=1}^n [ {\ele}, \dots, \overbrace{{\eme}}^i, \dots,  {\ele} ]} \to 0
\]}
\end{Pro}
{\it Proof.}
Consider the following commutative diagram of Hom-Leibniz $n$-algebras where $\pi$ denotes the canonical projection on the quotient
\[ \xymatrix{
0 \ar[d] & 0 \ar[d]\\
( \displaystyle \bigoplus_{i=1}^n {\ele} \ast \dots \ast \overbrace{{\eme}}^i \ast \dots \ast {\ele},\widetilde{{\alpha}}_{{\ele} \ast \dots \ast {\ele}_{\mid}}) \ar[d] \ar[r]^{\quad  \quad \quad  \quad \quad  \quad {\psi}_{\mid}}& ({\eme}, \widetilde{\alpha}_{\eme}) \ar[d] \\
({\ele} \ast \dots \ast {\ele},\widetilde{{\alpha}}_{{\ele} \ast \dots \ast {\ele}}) \ar[d]^{\pi \ast \dots \ast \pi} \ar[r]^{\quad \quad \quad \quad  \quad \quad  \psi} & ({\ele},\widetilde{\alpha}_{\ele}) \ar[d]^{\pi}\\
(\frac{\ele}{\eme} \ast \dots \ast \frac{\ele}{\eme},\widetilde{\alpha}_{{\ele}/{\eme} \ast \dots \ast {\ele}/{\eme}}) \ar[d] \ar[r]^{\quad \quad \quad  \quad \quad  \quad  \overline{\psi}}& (\frac{\ele}{\eme},\overline{\widetilde{\alpha}}_{{\ele}})  \ar[d]\\
0 & 0
} \]
where ${\psi} (l_1 \ast \dots \ast l_n)=[l_1, \dots,l_n]$. Then, forgetting the Hom-Leibniz $n$-algebra structures, by using the Snake Lemma for the same diagram of vector spaces, we obtain the following exact sequence,
$$\Ker(\psi_{\mid}) \to \Ker(\psi) \to \Ker(\overline{\psi}) \to \cok(\psi_{\mid}) \to \cok(\psi) \to \cok(\overline{\psi}) \to 0$$
where  $\Ker\left(  \psi_{\mid}\right)  \cong \Ker( \displaystyle \bigoplus_{i=1}^n {\ele} \ast \dots \ast \overbrace{{\eme}}^i \ast \dots \ast {\ele}\to {\eme})$; $\Ker\left(  \psi\right)  \cong {_n}HL_1^{\alpha}({\ele})$ and $\Ker\left(  \overline{\psi}\right)  \cong {_n}HL_1^{\alpha}({\ele}/{\eme})$ by Remark \ref{nucleo ext central Leib}; $\cok\left(  \psi_{\mid}\right)  \cong \frac{{\eme}}{\displaystyle \bigoplus_{i=1}^n[{\ele}, \dots , \overbrace{{\eme}}^i, \dots, {\ele} ]}$ and $\cok\left(  \psi\right)  = \cok\left(  \overline{\psi}\right)  =0$. \rdg


\section{Unicentrality of Hom-Leibniz $n$-algebras}

Our goal in this section is the generalization of the concept and properties of unicentral Leibniz algebras to the setting of Hom-Leibniz $n$-algebras. Namely  (see \cite{CC}), a Leibniz algebra $\frak{q}$ is said to be unicentral if  $\pi(Z({\frak{g}}) )= Z(\frak{q})$ for every central extension $\pi : \frak{g} \twoheadrightarrow \frak{q}$. In particular, perfect Leibniz algebras are unicentral (see \cite[Proposition 4]{CC}).

As a first step, we show in the following example, that central extensions of perfect Hom-Leibniz $n$-algebras are not unicentral.

\begin{Ex} \label{unicentral}
Let $({\ele}, \widetilde{\alpha}_{\ele})$ be the three-dimensional Hom-Leibniz $3$-algebra with basis $\{e_1, e_2, e_3 \}$, bracket operation given by $[e_1, e_1,e_1]=e_1; [e_1, e_1, e_2]=e_2;$ $[e_1,e_2,e_1]=e_3$  and zero elsewhere, and $\widetilde{\alpha}_{\ele} = (0)$. Obviously, $({\ele}, \widetilde{\alpha}_{\ele})$ is a perfect Hom-Leibniz $3$-algebra.

Consider the four-dimensional Hom-Leibniz $3$-algebra  $({\ka}, \widetilde{\alpha}_{\ka})$ with basis $\{a_1,$ $a_2,a_3, a_4 \}$, bracket operation given by $[a_3, a_3,a_3]=a_3; [a_3, a_3, a_1]=a_1; [a_3, a_1, a_3] = a_2; [a_3, a_3, a_2] = a_4$  and zero elsewhere, and $\widetilde{\alpha}_{\ka} = (0)$.

The surjective homomorphism $f : ({\ka}, \widetilde{\alpha}_{\ka}) \twoheadrightarrow ({\ele}, \widetilde{\alpha}_{\ele})$ given by $f(a_1)=e_2; f(a_2)$ $= e_3; f(a_3) = e_1; f(a_4) = 0$, is a central extension since $\Ker(f) = \langle \{a_4 \} \rangle$ and $Z({\ka}) = \langle \{ a_4 \} \rangle$. Moreover, $f(Z({\ka})) = 0$, but $Z({\cal L}) = \langle \{e_3 \}\rangle$, hence $f(Z({\ka})) \varsubsetneqq Z({\cal L})$.
\end{Ex}

By this fact, in what follows we show some results concerning the generalization of properties of unicentral Leibniz algebras.

\begin{De}
A perfect Hom-Leibniz $n$-algebra $({\ele}, \widetilde{\alpha}_{\ele})$ is said to be centrally closed if its universal central extension is
$$0 \to 0 \to ({\ele}, \widetilde{\alpha}_{\ele}) \stackrel{\sim}\to ({\ele}, \widetilde{\alpha}_{\ele}) \to 0$$
i.e. ${_n}HL_1^{\alpha}({\ele}) = 0$ and $(\frak{uce}({\ele}), \widetilde{\alpha}_{\frak{uce}({\ele})}) \cong ({\ele}, \widetilde{\alpha}_{\ele})$
\end{De}

\begin{Le} \label{lema 5}
Let $f : ({\ka}, \widetilde{\alpha}_{\ka}) \twoheadrightarrow ({\ele}, \widetilde{\alpha}_{\ele})$ be a central extension of a perfect Hom-Leibniz $n$-algebra $({\ele}, \widetilde{\alpha}_{\ele})$. Then the following statements hold:
\begin{enumerate}
\item[a)] ${\ka} = [{\ka}, \dots, {\ka}] + \Ker(f)$.
\item[b)] If $\alpha_{\ele}(l) \in Z(\alpha_{\ele}({\ele}))$, then $[l_1, \dots, l_{i-1},l,l_{i+1},\dots,l_n] \in \Ker(\alpha_{\ele})$, for all $l_j \in {\ele}, i \in \{1, 2, \dots, n \}, j \in \{1, \dots, \widehat{i}, \dots, n \}$.
\end{enumerate}
\end{Le}
{\it Proof.} {\it a)} For any $k \in {\ka}$, $f(k) \in {\ele} = [{\ele}, \dots, {\ele}]$, then $f(k) = [f(k_1), \dots, f(k_n)]$, hence $k - [k_1, \dots, k_n] \in \Ker(f)$.

{\it b)} If $\alpha_{\ele}(l) \in Z(\alpha_{\ele}({\ele}))$, then $[\alpha_{\ele}(l_1), \dots, \alpha_{\ele}(l),\dots, \alpha_{\ele}(l_n)]=0$.  \rdg

\begin{Pro} \label{prop 1}
Let $f : ({\ka}, \widetilde{\alpha}_{\ka}) \twoheadrightarrow ({\ele}, \widetilde{\alpha}_{\ele})$ be a central extension of a perfect Hom-Leibniz $n$-algebra $({\ele}, \widetilde{\alpha}_{\ele})$ with $\alpha_{\ele}$ injective, such that $({\ka}, \widetilde{\alpha}_{\ka})$ satisfies the following condition
\begin{equation} \label{condition}
[\alpha(k), \alpha(k), \alpha(k_3), \dots, \alpha(k_n)]=0, \text{for all}\ k, k_3, \dots, k_n \in {\ka}
\end{equation}
Then $$f(Z(\alpha_{\ka}({\ka}))) = Z(\alpha_{\ele}({\ele}))$$
\end{Pro}
{\it Proof.} Let $\alpha_{\ka}(k) \in Z(\alpha_{\ka}({\ka}))$, then $f(\alpha_{\ka}(k)) \in Z(\alpha_{\ele}({\ele}))$ since
\[
\begin{array}{l}
\left[\alpha_{\ele}(l_1), \dots, f(\alpha_{\ka}(k)), \dots , \alpha_{\ele}(l_n) \right] = \\

 [\alpha_{\ele}(f(k_1)),\dots ,f(\alpha_{\ka}(k)), \dots, \alpha_{\ele}(f(k_n))] = \\

 f[\alpha_{\ka}(k_1), \dots, \alpha_{\ka}(k), \dots, \alpha_{\ka}(k_n)] = 0
 \end{array}
\]

Conversely, for any $\alpha_{\ele}(l) \in Z(\alpha_{\ele}({\ele}))$, there exists any $k \in {\ka}$ such that $f(k)=l$, hence $\alpha_{\ele}(l) = \alpha_{\ele}(f(k)) = f(\alpha_{\ka}(k))$. We must show that $\alpha_{\ka}(k) \in Z(\alpha_{\ka}({\ka}))$. Indeed,
\[
\begin{array}{l}
[\alpha_{\ka}(k), \alpha_{\ka}(k_2), \dots, \alpha_{\ka}(k_n)] =\\
 - [\alpha_{\ka}(k_2), \alpha_{\ka}(k), \alpha_{\ka}(k_3),\dots, \alpha_{\ka}(k_n)] = \\
- [\alpha_{\ka}[k_{21}, \dots, k_{2n}]+ \Ker(f), \alpha_{\ka}(k), \alpha_{\ka}(k_3), \dots, \alpha_{\ka}(k_{n})]
\end{array}
\]
by condition (\ref{condition}) and  Lemma \ref{lema 5} {\it a)}. Applying the fundamental identity (\ref{def}) and having in mind that $\Ker(f) \subseteq Z({\ka})$, the above equality reduces to
\[
\begin{array}{l}
- [[\alpha_{\ka}(k_{21}), k, k_3, \dots, k_n],\alpha_{\ka}^2(k_{22}), \dots, \alpha_{\ka}^2(k_{2n})]\\

- [\alpha_{\ka}^2(k_{21}), [\alpha_{\ka}(k_{22}), k, k_3, \dots, k_n], \alpha_{\ka}^2(k_{23}), \dots, \alpha_{\ka}^2(k_{2n})] - \dots \\

- [\alpha_{\ka}^2(k_{21}), \dots, \alpha_{\ka}^2(k_{2(n-1)}), [\alpha_{\ka}(k_{2n}), k, k_3, \dots, k_n]]\\
\end{array}
\]
which vanishes since the brackets of the form $[\alpha_{\ka}(k_{2i}), k, k_3, \dots, k_n]$ are in $\Ker(f) \subseteq Z({\ka})$ because $f[\alpha_{\ka}(k_{2i}), k, k_3, \dots, k_n] = [f(\alpha_{\ka}(k_{2i})), l, f(k_3), \dots, f(k_n)] \in \Ker(\alpha_{\ele})$ by Lemma \ref{lema 5} {\it b)}, and $\alpha_{\ele}$ is injective.

The vanishing of the  other possible brackets is completely analogous to the last arguments, so we omit it. \rdg

\begin{Rem}
Hom-Lie $n$-algebras are examples of Hom-Leibniz $n$-algebras satisfying condition (\ref{condition}). Also Hom-Lie triple systems satisfy condition (\ref{condition}) in case $n=3$ (see Example \ref{ejemplo 1} (d)).
\end{Rem}

\begin{Ex} In the following we present an example of central extension satisfying the conditions established in Proposition \ref{prop 1}.

Consider the four-dimensional $\mathbb{C}$-vector space ${\ele}$ with basis $\{e_1,e_2,$ $e_3, e_4\}$ endowed with the ternary bracket operation given by $[e_2,e_3,e_4]=e_1;$ $[e_1,e_3,e_4]=e_2; [e_1,e_2,e_4]=e_3;[e_1,e_2,e_3]=e_4$, together with the corresponding skew-symmetric ones and zero elsewhere. By Lemma 2.2 in \cite{BZWL}, $({\ele}, [-,-,-])$ is a Lie 3-algebra.

Now consider the homomorphism of Lie 3-algebras $\alpha : {\ele} \to {\ele}$ given by $\alpha(e_1)=e_1; \alpha(e_2)= - e_2; \alpha(e_3)=e_3; \alpha(e_4)= - e_4$. Then Theorem 3.4 in \cite{AMS} endows ${\ele}$ with a structure of Hom-Leibniz 3-algebra with bracket operation given by $\{e_2,e_3,e_4\}=  e_1; \{e_1,e_3,e_4\}= - e_2; \{e_1,e_2,e_4\}=  e_3;\{e_1,e_2,e_3\}= - e_4$, together with the corresponding skew-symmetric ones and zero elsewhere, and $\widetilde{\alpha}_{\ele} = (\alpha, \alpha)$. This Hom-Leibniz 3-algebra is perfect and $\alpha_{\ele}$ is injective.

Consider the four-dimensional $\mathbb{C}$-vector space ${\ka}$ with basis $\{a_1,a_2,a_3, a_4\}$ endowed with the ternary bracket operation given by $[a_2,a_3,a_4]=a_1;[a_1,a_3,a_4]=a_2; [a_1,a_2,a_4]=a_3;[a_1,a_2,a_3]=a_4$, together with the corresponding skew-symmetric ones and zero elsewhere. By Lemma 2.2 in \cite{BZWL}, $({\ka}, [-,-,-])$ is a Lie 3-algebra.

Now consider the homomorphism of Lie 3-algebras $\beta : {\ka} \to {\ka}$ given by $\beta(a_1)=- a_1; \beta(a_2)= a_2; \beta(a_3)=- a_3; \beta(a_4)= a_4$. Then Theorem 3.4 in \cite{AMS} endows ${\ka}$ with a structure of Hom-Leibniz 3-algebra with bracket operation given by $\{a_2,a_3,a_4\}=  - a_1; \{a_1,a_3,a_4\}= a_2; \{a_1,a_2,a_4\}=  - a_3;\{a_1,a_2,a_3\}= a_4$, together with the corresponding skew-symmetric ones and zero elsewhere, and $\widetilde{\alpha}_{\ka} = (\beta, \beta)$. This Hom-Leibniz 3-algebra obviously satisfies condition (\ref{condition}).

The surjective homomorphism $f : ({\ka}, \widetilde{\alpha}_{\ka}) \twoheadrightarrow ({\ele}, \widetilde{\alpha}_{\ele})$ defined by $f(a_1) = e_2; f(a_2) = e_1; f(a_3) =  e_4; f(a_4) = - e_3$, is a central extension since $\Ker(f)$ and $Z({\ka})$ are both trivial.
\end{Ex}

Let $({\ele}, \widetilde{\alpha}_{\ele})$ be a perfect Hom-Leibniz $n$-algebra with $\alpha_{\ele}$ injective. Assume that $({\ele}, \widetilde{\alpha}_{\ele})$ satisfies condition (\ref{condition}). Then $\left( \frak{uce}({\ele}), \widetilde{\alpha}_{\frak{uce}(\ele)} \right)$ satisfies condition (\ref{condition}) provided that $\{ \alpha_{\ele}(l), \alpha_{\ele}(l), \alpha_{\ele}(l_3), \dots, \alpha_{\ele}(l_n) \} \in {_n}HL_1^{\alpha}({\ele})$, for all $l, l_3, \dots, l_n \in {\ele}$, is the zero coset. This fact occurs, for instance, when  $({\ele}, \widetilde{\alpha}_{\ele})$ is centrally closed.

From now on, we assume that $\left( \frak{uce}({\ele}), \widetilde{\alpha}_{\frak{uce}(\ele)} \right)$ satisfies condition (\ref{condition}) when $({\ele}, \widetilde{\alpha}_{\ele})$ does. Then Proposition \ref{prop 1} gives the following equality:
\begin{equation} \label{eq}
u_{\ele}(Z( \alpha_{\frak{uce({\ele})}} (\frak{uce}({\ele})))) = Z(\alpha_{\ele}({\ele}))
\end{equation}

\begin{Th} \label{isomorph}
Let $({\ele}, \widetilde{\alpha}_{\ele})$ be a perfect Hom-Leibniz $n$-algebra with $\alpha_{\ele}$ injective such that  $({\ele}, \widetilde{\alpha}_{\ele})$ and  $\left( \frak{uce}({\ele}), \widetilde{\alpha}_{\frak{uce({\ele})}} \right)$  satisfy condition (\ref{condition}). Then there is an isomorphism
$$\frac{\alpha_{\ele}({\ele})}{Z( \alpha_{\ele}({\ele}))} \cong \frac{\alpha_{\frak{uce({\ele})}}(\frak{uce}({\ele}))}{Z(\alpha_{\frak{uce({\ele})}}({\frak{uce}(\ele)}))}$$
\end{Th}
{\it Proof.} The universal central extension of $({\ele}, \widetilde{\alpha}_{\ele})$ induces the central extension
$$0\to \left(\alpha_{\frak{uce}({\ele})} \left( _{n}HL_{1}^{\alpha }\left( {\ele} \right) \right), \widetilde{\alpha}_{\frak{uce}(\ele)\mid}\right) \to \left( \alpha_{\frak{uce({\ele})}} (\frak{uce}({\ele})), \widetilde{\alpha}_{\frak{uce}(\ele) \mid} \right) \to  \left( \alpha_{\ele}({\ele}),\widetilde{\alpha }_{\ele \mid}\right) \to 0$$
when $\alpha_{\ele}$ is injective. Moreover condition (\ref{condition}) is preserved by the terms in this central extension.

Bearing in mind (\ref{eq}), the kernels of the horizontal arrows in the commutative diagram
\[
\xymatrix{
Z(\alpha_{\frak{uce}({\ele})}({\frak{uce}(\ele)})) \ar@{>>}[r]^{\quad u_{{\ele}\mid}} \ar@{>->}[d] & Z( \alpha_{\ele}({\ele})) \ar@{>->}[d]\\
 \alpha_{\frak{uce}({\ele})} (\frak{uce}({\ele})) \ar@{>>}[r]^{\quad u_{{\ele}\mid}} &   \alpha_{\ele}({\ele})
} \]
coincide, then the cokernels of the vertical homomorphisms are isomorphic. \rdg

\begin{Pro} \label{composition}
Let $0 \to ({\eme}, \widetilde{\alpha}_{\eme}) \to ({\ka}, \widetilde{\alpha}_{\ka}) \stackrel{\pi} \to ({\ele}, \widetilde{\alpha}_{\ele}) \to 0$ and $0 \to ({\ene}, \widetilde{\alpha}_{\ene}) \to ({\he}, \widetilde{\alpha}_{\he}) \stackrel{\tau} \to ({\ka}, \widetilde{\alpha}_{\ka}) \to 0$ be central extensions of Hom-Leibniz $n$-algebras. Then the following statements hold:
\begin{enumerate}
\item[a)] If $\pi \circ \tau : ({\he}, \widetilde{\alpha}_{\he}) \twoheadrightarrow ({\ele}, \widetilde{\alpha}_{\ele})$ is a universal $\alpha$-central extension, then $\tau : ({\he}, \widetilde{\alpha}_{\he}) \twoheadrightarrow ({\ka}, \widetilde{\alpha}_{\ka})$ is a universal central extension.
\item[b)] If $\tau : ({\he}, \widetilde{\alpha}_{\he}) \twoheadrightarrow ({\ka}, \widetilde{\alpha}_{\ka})$ is a universal central extension, then  $\pi \circ \tau : ({\he}, \widetilde{\alpha}_{\he}) \twoheadrightarrow ({\ele}, \widetilde{\alpha}_{\ele})$ is an $\alpha$-central extension which is universal over central extensions, that is, for any central extension $0 \to ({\as}, \widetilde{\alpha}_{\as}) \to ({\pe}, \widetilde{\alpha}_{\pe}) \stackrel{\omega} \to ({\ele}, \widetilde{\alpha}_{\ele}) \to 0$ there exists a unique homomorphism $\Phi : ({\he}, \widetilde{\alpha}_{\he}) \to ({\pe}, \widetilde{\alpha}_{\pe})$ such that $\omega \circ \Phi = \pi \circ \tau$.
\end{enumerate}
\end{Pro}
{\it Proof.} {\it a)} If  $\pi \circ \tau : ({\he}, \widetilde{\alpha}_{\he}) \twoheadrightarrow ({\ele}, \widetilde{\alpha}_{\ele})$ is a universal $\alpha$-central extension, then ${_n}HL_0^{\alpha}({\he}) = {_n}HL_1^{\alpha}({\he}) = 0$ by Theorem \ref{teorema}
 {\it d)}. Hence  $\tau : ({\he}, \widetilde{\alpha}_{\he}) \twoheadrightarrow ({\ka}, \widetilde{\alpha}_{\ka})$ is a universal central extension by Theorem \ref{teorema} {\it e)}.

 {\it b)} If $\tau : ({\he}, \widetilde{\alpha}_{\he}) \twoheadrightarrow ({\ka}, \widetilde{\alpha}_{\ka})$ is a universal central extension, then $({\he}, \widetilde{\alpha}_{\he})$ and $({\ka}, \widetilde{\alpha}_{\ka})$ are perfect Hom-Leibniz $n$-algebras by Lemma \ref{lema 3}, hence $\pi \circ \tau : ({\he}, \widetilde{\alpha}_{\he}) \twoheadrightarrow ({\ele}, \widetilde{\alpha}_{\ele})$ is an $\alpha$-central extension by Lemma \ref{lema 4}. Moreover  $\pi \circ \tau$ is universal over central extensions. Indeed, for any central extension $0 \to ({\as}, \widetilde{\alpha}_{\as}) \to ({\pe}, \widetilde{\alpha}_{\pe}) \stackrel{\omega} \to ({\ele}, \widetilde{\alpha}_{\ele}) \to 0$, construct the pull-back extension
 \[ \xymatrix{
  0 \ar[r] & ({\as}, \widetilde{\alpha}_{\as}) \ar[r] \ar@{=}[d]&  ({\pe} \times_{\ele} {\ka}, \widetilde{\alpha}_{\pe}\times \widetilde{\alpha}_{\ka}) \ar[r]^{\quad \quad \overline{\omega}} \ar[d]^{\overline{\pi}}& ({\ka}, \widetilde{\alpha}_{\ka}) \ar[r] \ar[d]^{\pi}& 0\\
 0 \ar[r] & ({\as}, \widetilde{\alpha}_{\as}) \ar[r] &  ({\pe}, \widetilde{\alpha}_{\pe}) \ar[r]^{\omega} & ({\ele}, \widetilde{\alpha}_{\ele}) \ar[r] & 0
 } \]
which is central. Since $\tau : ({\he}, \widetilde{\alpha}_{\he}) \twoheadrightarrow ({\ka}, \widetilde{\alpha}_{\ka})$ is a universal central extension, then there exists a unique homomorphism $\varphi :  ({\he}, \widetilde{\alpha}_{\he}) \to ({\pe} \times_{\ele} {\ka}, \widetilde{\alpha}_{\pe}\times \widetilde{\alpha}_{\ka})$ such that $\overline{\omega} \circ \varphi = \tau$. Then $\Phi = \overline{\pi} \circ \varphi$ satisfies the required universal property thanks to Lemma \ref{lema 2}. \rdg

\begin{Co}
Let $({\ele}, \widetilde{\alpha}_{\ele}), ({\ele}', \widetilde{\alpha}_{\ele'})$ be  perfect Hom-Leibniz $n$-algebras with both $\alpha_{\ele}, \alpha_{\ele'}$ injective and such that  $({\ele}, \widetilde{\alpha}_{\ele}), ({\ele}', \widetilde{\alpha}_{\ele'}), \left( \frak{uce}({\ele}), \widetilde{\alpha}_{\frak{uce}({\ele})} \right)$, and $\left( \frak{uce}({\ele}'), \right.$ $\left.\widetilde{\alpha}_{\frak{uce}({\ele}')} \right)$  satisfy condition (\ref{condition}). Then
 \begin{enumerate}
 \item[a)] If  $\left( \frak{uce}({\ele}), \widetilde{\alpha}_{\frak{uce}({\ele})} \right) \cong \left( \frak{uce}({\ele}'), \widetilde{\alpha}_{\frak{uce}({\ele}')} \right)$, then $\frac{\alpha_{\ele}({\ele})}{Z( \alpha_{\ele}({\ele}))} \cong \frac{\alpha_{\ele'}({\ele'})}{Z( \alpha_{\ele'}({\ele'}))}$.
 \item[b)] If $\frac{\alpha_{\ele}({\ele})}{Z( \alpha_{\ele}({\ele}))} \cong \frac{\alpha_{\ele'}({\ele'})}{Z( \alpha_{\ele'}({\ele'}))}$, then $\left( \frak{uce}(\alpha_{\ele}({\ele})), \widetilde{\alpha}_{\frak{uce}({\ele})\mid} \right) \cong \left( \frak{uce}(\alpha_{\ele'}({\ele'})), \widetilde{\alpha}_{\frak{uce}({\ele}')\mid} \right)$.
 \end{enumerate}
\end{Co}
{\it Proof.} {\it a)}
 If  $\left( \frak{uce}({\ele}), \widetilde{\alpha}_{\frak{uce}({\ele})} \right) \cong \left( \frak{uce}({\ele}'), \widetilde{\alpha}_{\frak{uce}({\ele}')} \right)$, then $\frac{\alpha_{\ele}({\ele})}{Z( \alpha_{\ele}({\ele}))} \cong \frac{\alpha_{\frak{uce}({\ele})}(\frak{uce}({\ele}))}{Z(\alpha_{\frak{uce}({\ele})}({\frak{uce}(\ele)}))} \cong \frac{\alpha_{\frak{uce}({\ele}')}(\frak{uce}({\ele}'))}{Z(\alpha_{\frak{uce}({\ele}')}({\frak{uce}(\ele')}))}$ $\cong \frac{\alpha_{\ele'}({\ele'})}{Z( \alpha_{\ele'}({\ele'}))}$ by Theorem \ref{isomorph}.

{\it b)} If $\frac{\alpha_{\ele}({\ele})}{Z( \alpha_{\ele}({\ele}))} \cong \frac{\alpha_{\ele'}({\ele'})}{Z( \alpha_{\ele'}({\ele'}))}$, then $\left( \frak{uce} \left( \frac{\alpha_{\ele}({\ele})}{Z( \alpha_{\ele}({\ele}))} \right), \widetilde{\alpha}_{\frak{uce}({\ele})} \right) \cong \left( \frak{uce} \left( \frac{\alpha_{\ele'}({\ele'})}{Z( \alpha_{\ele'}({\ele'}))}\right),\right.$ $\left.\widetilde{\alpha}_{\frak{uce}({\ele}')}  \right)$.

\noindent Now, applying Proposition \ref{composition} {\it b)} to the central extensions $u_{\ele \mid} : \left( \frak{uce}(\alpha_{\ele}({\ele})), \widetilde{\alpha}_{\frak{uce}{\ele})\mid} \right)$ $\to (\alpha_{\ele}({\ele}), \widetilde{\alpha}_{\ele \mid})$ and $ p : (\alpha_{\ele}({\ele}), \widetilde{\alpha}_{\ele \mid}) \to (\alpha_{\ele}({\ele})/Z(\alpha_{\ele}({\ele})), \widetilde{\overline{\alpha}}_{\ele})$, we conclude that $\left( \frak{uce}(\alpha_{\ele}({\ele})), \widetilde{\alpha}_{\frak{uce}{\ele})\mid} \right) \cong \left( \frak{uce}(\alpha_{\ele}({\ele})/Z(\alpha_{\ele}({\ele}))), \widetilde{\overline{\alpha}}_{\frak{uce}({\ele})\mid} \right)$. \rdg

\begin{Co}
Let $({\ele}, \widetilde{\alpha}_{\ele})$ be a centerless perfect Hom-Leibniz $n$-algebra with $\alpha_{\ele}$ injective such that  $({\ele}, \widetilde{\alpha}_{\ele})$ and  $\left( \frak{uce}({\ele}), \widetilde{\alpha}_{\frak{uce}({\ele})} \right)$  satisfy condition (\ref{condition}). Then $Z(\alpha_{\frak{uce}({\ele})}(\frak{uce}({\ele}))) \cong {_n}HL_1^{\alpha}(\alpha_{\ele}({\ele}))$ and the universal central extension of  $(\alpha_{\ele}({\ele}),$ $\widetilde{\alpha}_{\ele \mid})$ is
$$0 \to \left( Z(\alpha_{\frak{uce}({\ele})}(\frak{uce}({\ele}))), \widetilde{\alpha}_{\frak{uce}{\ele})\mid} \right) \to  \left( \frak{uce}(\alpha_{\ele}({\ele})), \widetilde{\alpha}_{\frak{uce}({\ele}) \mid} \right) \to (\alpha_{\ele}({\ele}), \widetilde{\alpha}_{\ele \mid}) \to 0$$
\end{Co}
{\it Proof.}  $Z({\ele})=0$ and $\alpha_{\ele}$ injective implies that $Z(\alpha_{\ele}({\ele}))=0$. Then Theorem \ref{isomorph} implies that $0 \to \left( Z(\alpha_{\frak{uce}({\ele})}(\frak{uce}({\ele}))), \widetilde{\alpha}_{\frak{uce}({\ele})\mid} \right) \to  \left( \alpha_{\ele}(\frak{uce}({\ele})), \widetilde{\alpha}_{\frak{uce}({\ele}) \mid} \right) \to (\alpha_{\ele}({\ele}), \widetilde{\alpha}_{\ele \mid}) \to 0$ is isomorphic to the universal central extension of $(\alpha_{\ele}({\ele}), \widetilde{\alpha}_{\ele \mid})$. \rdg


\section*{Acknowledgements}

First author was supported by Ministerio de Economía y Competitividad (Spain) (European FEDER support included), grant MTM2013-43687-P.


\begin{center}

\end{center}

\end{document}